\documentclass[11pt]{amsart}
\usepackage{Macros}

\title[Refined Cyclic Sieving on Words for the Major Index Statistic]
  {Refined Cyclic Sieving on Words \\ for the Major Index Statistic}
\author{Connor Ahlbach and Joshua P. Swanson}
\address{Department of Mathematics, University of Washington,
Seattle, WA 98195, USA}
\thanks{Accepted to the European Journal of Combinatorics. \copyright 2018. This manuscript version is made available under the CC-BY-NC-ND 4.0 license \url{http://creativecommons.org/licenses/by-nc-nd/4.0/}.}
\date{\today}

\begin{document}

\begin{abstract}
  Reiner--Stanton--White \cite{MR2087303} defined the cyclic sieving phenomenon (CSP) associated to a finite cyclic group action and a polynomial. A key example arises from the length generating function for minimal length coset representatives of a parabolic quotient of a finite Coxeter group. In type $A$, this result can be phrased in terms of the natural cyclic action on words of fixed content.

  There is a natural notion of \textit{refinement} for many CSP's. We formulate and prove a refinement, with respect to the major index statistic, of this CSP on words of fixed content by also fixing the \textit{cyclic descent type}. The argument presented is completely different from Reiner--Stanton--White's representation-theoretic approach. It is combinatorial and largely, though not entirely, bijective in a sense we make precise with a ``universal'' sieving statistic on words, $\flex$.

  A building block of our argument involves cyclic sieving for shifted subset sums, which also appeared in Reiner--Stanton--White. We give an alternate, largely bijective proof of a refinement of this result by extending some ideas of Wagon--Wilf \cite{MR1269164}.
\end{abstract}
\maketitle

\section{Introduction}
\label{sec:intro}

Since Reiner, Stanton, and White introduced the \textit{cyclic sieving phenomenon} (CSP) in 2004 \cite{MR2087303}, it has become an important companion to any cyclic action on a finite set. Some remarkable examples of the CSP involve the action of a Springer regular element on Coxeter groups \cite[Theorem~1.6]{MR2087303}, the action of Schutzenberger's promotion on Young tableaux of fixed rectangular shape \cite{MR2557880}, and the creation of new CSPs from old using multisets and plethysms with homogeneous symmetric functions \cite[Proposition 8]{MR2837599}. See \cite{MR2866734} for Sagan's thorough introduction to the cyclic sieving phenomenon. More recent work on the CSP includes \cite{MR3281144,MR3207480,MR3537922}. Here we are concerned with cyclic sieving phenomena involving cyclic descents on words. Cyclic descents were used implicitly by Klyachko \cite{klyachko74} and independently introduced by Cellini \cite{MR1637728}. Since then, cyclic descents have been used by Lam and Postnikov in studying alcoved polytopes \cite{1202.4015} and by Petersen in studying $P$-partitions \cite{MR2181371}. They also appear prominently in an ongoing line of research on cyclic descent extensions for standard tableaux by Adin, Elizalde, Reiner, and Roichman \cite{MR3682729,1710.06664,1801.00044}.

An earlier ``extended abstract'' for the present work appeared in \cite{as17}. We assume some familiarity with the CSP, though we recall certain key statements.

\begin{defn}
  \label{def:csp1}
  Suppose $C_n$ is a cyclic group of order $n$ generated by $\sigma_n$, $W$ is a finite set on which $C_n$ acts, and $f(q) \in \bN[q]$. We say the triple
  $(W, C_n, f(q))$ exhibits the
  \textit{cyclic sieving phenomenon (CSP)} if for all $k \in \bZ$,
  \begin{align}
    \label{CSP1}
      \# W^{\si_n^k} \coloneqq \# \{ w \in W: \si_n^k \dd w = w \} = f(\w_n^k),
  \end{align}
  where $ \w_n $ is any fixed primitive $ n $-th root of unity.
\end{defn}
\noindent Representation theoretically, evaluations of $f$ at $ n $-th roots of unity yield the characters of the $C_n$-action on $W$.

In many instances of cyclic sieving, and all of those considered here, $f(q)$ is the generating function for some statistic on $ W $. Given a statistic $\stat \colon W \to \bN$,
let
\begin{align}
\label{eqn:statGF}
    W^{\stat}(q) \coloneqq \sum_{w \in W} q^{\stat w} \in \bN[q].
\end{align}
We say two statistics $ \stat, \stat' \colon W \to \bN $ are \emph{equidistributed} on $ W $ if $ W^{\stat}(q) = W^{\stat'}(q) $, and we say they are \emph{equidistributed modulo $ n $} on $ W $ if $ W^{\stat}(q) \equiv W^{\stat'}(q) \; (\Mod q^n - 1) $.

Our main result is a refinement of a CSP triple first observed by
Reiner--Stanton--White, which we now summarize; see Section~\ref{sec:comb} for missing definitions. Consider words in the alphabet $\bP \coloneqq \{1, 2, \ldots\}$. Given a word $w = w_1\cdots w_n$ of length $n$, let $\cont(w)$ denote the \textit{content} of $w$ and write
\begin{align}
\label{eqn:FixedCont}
      \W_\alpha \coloneqq \{\text{words }w : \cont(w) = \alpha\}
\end{align}
for the set of words with content $\alpha$. Write $\maj(w)$ for the \textit{major index} of $w$. The cyclic group $C_n$ acts on words of length $n$ by rotation.  The following expresses an interesting result of Reiner, Stanton, and White in our notation.

\begin{thm}{\cite[Proposition~4.4]{MR2087303}}.
  \label{thm:rsw_alpha}
  Let $\alpha \vDash n$. The triple
    \[ (\W_\alpha, C_n, \W_\alpha^{\maj}(q)) \]
  exhibits the CSP.
\end{thm}

Reiner, Stanton, and White deduced Theorem~\ref{thm:rsw_alpha} from the following more general result about Coxeter systems.

\begin{thm}{\cite[Theorem~1.6]{MR2087303}}.
  \label{thm:rsw_Coxeter}
  Let $ (W,S) $ be a finite Coxeter system and $ J \subseteq S $. Let
  $W_J$ be the corresponding parabolic subgroup, $W^J$ the set of minimal
  length representatives for left cosets $X \coloneqq W/W_J$, and
  $ X^\ell(q) \coloneqq \sum_{ w \in W^J} q^{\ell(w)} $. Let $ C $ be a cyclic
  subgroup of $W$ generated by a Springer regular element. Then
  $(X, C, X^\ell(q))$ exhibits the cyclic sieving phenomenon.
\end{thm}

Theorem~\ref{thm:rsw_alpha} follows from Theorem~\ref{thm:rsw_Coxeter} when $ W = S_n $ by identifying $ W/W_J $ with words of fixed content $ \al $, where $ \al $ is the
composition recording the lengths of consecutive subsequences of $ J $,
and $C$ is generated by an $n$-cycle. One must
also use the classical result of MacMahon that $\maj$ is equidistributed
with the inversion statistic on words, from which it follows that
$\W_{\al}^{\maj}(q) = X^\ell(q)$ \cite[Art.~6]{MR1506186}.

\begin{defn}
  A \textit{refinement} of a CSP triple $(W, C_n, W^{\stat}(q))$ is a CSP
  triple
    \[  (V, C_n, V^{\stat}(q)) \]
  where $V \subset W$ has the restricted $C_n$-action.
\end{defn}

If $(V, C_n, V^{\stat}(q))$ refines $(W, C_n, W^{\stat}(q))$, then so does $(U, C_n, U^{\stat}(q)) $ where $ U \coloneqq W - V $. Thus, a CSP refinement partitions $ W $ into smaller CSPs with the same statistic. If $W$ is an orbit, its only refinements are
$W$ and $\varnothing$. In Section~\ref{sec:flex}, we define a statistic on words,
\textit{flex}, which is universal in the sense that it refines to all $C_n$-orbits. Such universal statistics are essentially equivalent to the choice of a total ordering for each orbit $\cO$ of $W$.

We partition words of fixed content into fixed \textit{cyclic descent type} ($ \CDT $). One
computes $\CDT(w)$ by building up $w$ by adding all $1$'s, $2$'s,
$\ldots$, and counting the number of \textit{cyclic descents} introduced at
each step. For precise details, see Definition~\ref{def:cdt} and Example~\ref{ex:cdt}. We
write the set of words with fixed content $\alpha \vDash n$ and cyclic
descent type $\delta$ as
\begin{align}
\label{eqn:FixedCDT}
    \W_{\alpha,\delta} \coloneqq \{ w \in \W_n : \cont(w) = \alpha, \CDT(w) = \delta\}.
\end{align}

Our main result is the following.

\begin{thm}
  \label{thm:alpha_delta}
  Let $\alpha \vDash n $ and $ \de $ be any composition. The triple
    \[ (\W_{\alpha,\delta}, C_n, \W_{\alpha,\delta}^{\maj}(q)) \]
  refines the CSP triple $(\W_\alpha, C_n, \W_\alpha^{\maj}(q))$.
\end{thm}

\noindent It is not at all clear how to modify Reiner--Stanton--White's
representation-theoretic approach to Theorem~\ref{thm:rsw_alpha} to give
Theorem~\ref{thm:alpha_delta}, since $\W_{\alpha, \delta}$ is not closed under
the $S_n$-action. Finding a representation-theoretic interpretation of
Theorem~\ref{thm:alpha_delta} would be quite interesting.

In the course of proving Theorem~\ref{thm:alpha_delta}, we derive an explicit
product formula for $\W_{\alpha,\delta}^{\maj}(q)$ mod $(q^n - 1)$
involving $q$-binomial coefficients, Theorem~\ref{thm:Maj GF Mod n}.
The formula results in a
$q$-identity similar to the Vandermonde convolution identity; see
Corollary~\ref{cor:vandermonde}. The argument involves constructing
$\W_{\alpha, \delta}$ algorithmically by recursively building a certain
tree.

The two-letter case of Theorem~\ref{thm:alpha_delta} can be rephrased as follows. Fix $n \in \bZ_{\geq 1}$ and $k, b \in \bZ_{\geq 0}$. Let $\SSS_{k,b}$
  denote the set of subsets $\Delta$ of
  $\bZ/n$ of size $k$ where $\#\{i \in \Delta : i+1 \not\in \Delta\} = b$.
  Define the statistic $\mbs \colon \SSS_{k,b} \to \bN$ by identifying $\bZ/n$ with $\{1, \ldots, n\}$ and setting $\mbs(\Delta) \coloneqq
  \sum_{i \in \Delta : i+1 \not\in \Delta} i$, which sums the maximum of the cyclic blocks of $\Delta$.

\begin{cor}
 The triple
    \[ (\SSS_{k,b}, C_n, \SSS_{k,b}^{\mbs}(q)) \]
  exhibits the CSP.
\end{cor}

\begin{ex} When $ n = 5, k = 3, b = 2 $,
\[
    \SSS_{k,b} = \{ \{ 1, 2, 4 \}, \{ 2, 3, 5 \}, \{ 3, 4, 1\}, \{ 4, 5, 2 \}, \{ 5, 1, 3 \} \},
\]
which have $\mbs$ statistic $6, 8, 5, 7, 4$, respectively, so
$ \SSS_{k,b}^{\mbs}(q) = q^4 + q^5 + q^6 + q^7 + q^8 $. We then have
$ \SSS_{k,b}^{\mbs}(\w_5) = 0, \SSS_{k,b}^{\mbs}(1) = 5$, in agreement with
\eqref{CSP1}.

\end{ex}

Theorem~8.3 in \cite{MR2087303} and hence Theorem~\ref{thm:rsw_alpha} builds on a representation-theoretic result due to Springer \cite[Proposition 4.5]{MR0354894}. Our argument is highly combinatorial, but it is not entirely
bijective. Finding an explicit bijection would be quite interesting.
See Section~\ref{sec:flex} for more details.

A key building block of our proof of Theorem~\ref{thm:alpha_delta} involves cyclic
sieving on multisubsets and subsets, which was also first stated in \cite{MR2087303}. We describe refinements of these results as well,
Theorem~\ref{thm:FixedsizesCSPMult} and Theorem~\ref{thm:Fixedgcds}, restricting to certain gcd requirements in the subset case. We present a completely different inductive proof of our subset refinement in the spirit of our proof of Theorem~\ref{thm:alpha_delta}.
Both our proof of Theorem~\ref{thm:Fixedgcds} and Theorem~\ref{thm:alpha_delta} use
an extension lemma, Lemma~\ref{lem:extend.csp}, which allows us to extend
CSPs from smaller cyclic groups to larger ones.

The rest of the paper is organized as follows. In Section~\ref{sec:comb}, we recall combinatorial background. In Section~\ref{modular_periodicity}, we introduce the concept of modular periodicity and prove our extension lemma, Lemma~\ref{lem:extend.csp}. In Section~\ref{sec:CDT}, we define cyclic descent type. In Section~\ref{sec:runs}, we decompose words with fixed content and cyclic descent type and prove a product formula for $ \W_{\al,\de}^{\maj}(q) $ modulo $ q^n - 1 $, Theorem~\ref{thm:Maj GF Mod n}. Section~\ref{sec:CSP} uses the results of Section~\ref{sec:runs} to prove our main result, Theorem~\ref{thm:alpha_delta}. Section~\ref{sec:subsets} refines cyclic sieving on multisubsets and subsets with respect to shifted sum statistics. In Section~\ref{sec:flex}, we introduce the flex statistic and use it to reinterpret Theorem~\ref{thm:alpha_delta}.

\section{Combinatorial Background}
\label{sec:comb}

In this section, we briefly recall or introduce combinatorial notions on words and fix our notation. We use the alphabet of positive integers $\bP \coloneqq \{1, 2, \ldots\}$ throughout unless otherwise noted. We also write $\#S$ or $ |S| $ for the cardinality of a set $S$. For $ \stat \colon W \to \bN $, recall the notation
\begin{align*}
    W^{\stat}(q) \coloneqq \sum_{w \in W} q^{\stat(w)}.
\end{align*}
A \emph{word} $w$ of length $ n $ is a sequence $ w = w_1 w_2 \cdots w_n $ of \textit{letters} $w_i \in \bP$. Let $ |w| $ denote the length of a word $ w $. Let $\W_n$ denote the set of all words of length $n$. The \textit{descent set} of $w$ is $\Des(w) \coloneqq \{1 \leq i < n : w_i > w_{i+1}\}$, and the number of descents is $\des(w) \coloneqq \# \Des(w)$. The \textit{major index} of $w$ is $\maj(w) \coloneqq \sum_{i \in \Des(w)} i$. The \textit{cyclic descent set} of $w$ is $\CDes(w) \coloneqq \{1 \leq i \leq n : w_i > w_{i+1}\}$, where now the subscripts are taken mod $n$, and we write $\cdes(w) \coloneqq \#\CDes(w)$ for the number of cyclic descents. Any position $1 \leq i \leq n$ that is not a cyclic descent is a \textit{cyclic weak ascent}. The \textit{inversion number} of $w$ is $\inv(w) \coloneqq \#\{(i, j) : 1 \leq i < j \leq n \text{ and }w_i > w_j\}$. We use lower dots between letters to indicate cyclic descents and upper dots to indicate cyclic weak ascents throughout the paper as in the following example.

\begin{ex}
  If $w = 155.3.155.3. = 1\udot5\udot531\udot5\udot53 $, then $|w| = 8$, $\Des(w) = \{ 3, 4, 7 \}$, $ \des(w) = 3 $, $\CDes(w) = \{3, 4, 7, 8\}$, $ \cdes(w) = 4 $, $\maj(w) = 14$, and $\inv(w) = 9$.
\end{ex}

A \emph{composition} or \emph{weak composition} of $n$ is a sequence $\al = (\al_1, \ldots, \al_m)$ of non-negative integers summing to $n$, typically denoted $\alpha \vDash n$. A composition is \emph{strong} if $\alpha_i > 0$ for all $i$. The \textit{content} of a word $ w $, denoted $ \cont(w) $, is the sequence $ \al $ whose $ j $-th part is the number of $ j $'s in $ w $. For $ w \in \W_n $, $ \cont(w) $ is a weak composition of $ n $. We write
\begin{align*}
  \W_\alpha &\coloneqq \{w \in \W_n : \cont(w) = \alpha \}.
\end{align*}
The cyclic group $ C_n \coloneqq \langle \si_n \rangle$ of order $ n $ acts on
$ \W_n $ by \textit{rotation} as
\begin{equation*}
  \si_n \cdot w_1 \cdots w_{n-1} w_n \coloneqq w_n w_1 \cdots w_{n-1}.
\end{equation*}
Typically we consider $\sigma_n$ to be the long cycle $(1 2 \cdots n) \in S_n$.

The set of all words in $\bP$ is a monoid under concatenation. A word is
\textit{primitive} if it is not a power of a smaller word. Any non-empty word
$w$ may be written uniquely as $w = v^f$ for $f \geq 1$ with $ v $ primitive.
We call $ |v| $ the \textit{period} of $w$, written $ \period(w) $, and
$f$ the \textit{frequency} of $w$, written $\freq(w)$. An orbit of $ \W_n $
under rotation is a \emph{necklace}, usually denoted $[w]$. We have
$\period(w) = \#[w]$ and $ \freq(w) \dd \period(w) = |w| $. Content,
primitivity, period, frequency, and $\cdes$ are all constant on necklaces.

\begin{ex}
  \label{ex:necklace}
  The necklace of $w = 15531553 = (1553)^2$ is
  \[
    [w] \coloneqq \{15531553, 55315531, 53155315, 31553155 \} \subset \W_{(2, 0, 2, 0, 4)}
           \subset \W_8,
  \]
  which has period $4$, frequency $2$, and $\cdes$ $4$.
\end{ex}

Reiner--Stanton--White gave equivalent conditions for a triple
$(W, C_n, f(q))$ to exhibit the CSP. In place of \eqref{CSP1} in
Definition~\ref{def:csp1}, we may instead require
\begin{align}
  \label{CSP2}
    f(q) \equiv \sum_{\text{orbits }\cO \subset W } \;
      \frac{q^n - 1}{q^{n/|\cO|} - 1} \quad (\Mod q^n - 1),
\end{align}
where the sum is over all orbits $ \cO $ under the action of $ C_n $ on $ W $. Note that for $ d \mid n $,
\begin{align*}
    \frac{q^n - 1}{q^d - 1} = \sum_{i=0}^{n/d - 1} q^{di}
     \not\equiv 0 \; (\Mod q^n - 1).
\end{align*}
This means every $ C_n $-action on a finite set $ W $ gives rise to a CSP $ (W, C_n, f(q)) $, where $ f(q) $ is the right hand side of \eqref{CSP2}. We refer the interested reader to \cite[Proposition 2.1]{MR2087303} for the proof of the equivalence of \eqref{CSP1} and \eqref{CSP2}.

\begin{rem}
\label{rem:relatedCSPs}
If $ (V, C_n, f(q)) $ exhibits the CSP, then so do both of the triples
$ (V, C_g, f(q)) $ and $ (V, C_n, f(q^{-1})) $ when $ g \mid n $ by
\eqref{CSP1}. In the latter case we have relaxed the
constraint $f(q) \in \bN[q]$ to $f(q) \in \bN[q, q^{-1}]$, which does
no harm since \eqref{CSP1} involves evaluations at roots of unity.
Further, if $ (V, C_n, f(q)) $ and $ (W, C_n, h(q)) $ exhibit
the CSP, then $ (V \coprod W,C_n, f(q) + g(q)) $ and $ (V \times W, C_n, f(q)h(q)) $ exhibit the CSP, where
$ C_n $ acts on $ V \times W $ by $ \tau \dd (v,w)
\coloneqq (\tau \dd v, \tau \dd w) $ \cite[Prop.~2.2]{MR2837599}.
\end{rem}

For a set $S$, write
\begin{align}
\label{eqn:SubMultisubNotation}
  \ch{S}{k} &\coloneqq \{\text{all $k$-element subsets of $S$}\}, \\
  \mch{S}{k} &\coloneqq \{\text{all $k$-element multisubsets of $S$}\}.
\end{align}
Let $\alpha = (\alpha_1, \ldots, \alpha_m) \vDash n$. We use the following standard $q$-analogues:
\begin{align*}
  [n]_q &\coloneqq 1 + q + \cdots + q^{n-1} = \frac{q^n - 1}{q - 1}, \\
  [n]_q! &\coloneqq [n]_q [n-1]_q \cdots [1]_q, \\
  \binom{n}{\alpha}_q
    &\coloneqq \frac{[n]_q!}{[\alpha_1]_q! \cdots [\alpha_m]_q!} \in \bN[q], \\
  \mch{n}{k}_q
    &\coloneqq \binom{n+k-1}{k}_q \coloneqq \binom{n+k-1}{k, n-1}_q.
\end{align*}

We write $[a, b] \coloneqq \{i \in \bZ : a \leq i \leq b\}$. Observe that the cyclic group $ C_n = \la \si_n \ra $ of order $ n $ acts on $ [0, n - 1] $ by $ \si_n(i) \coloneqq i + 1 \, (\Mod n) $. This induces actions of $ C_n $ on $ \ch{[0,n - 1]}{k} $ and $ \mch{[0,n - 1]}{k} $ by acting on values in each subset or multisubset. For example, $ \si_4 \cdot \{0,0,0,2,2,3\} = \{0,1,1,1,3,3\} $. These actions, in slightly more generality, appear in one of the original, foundational CSP results as follows.

\begin{thm}{\cite[Thm.~1.1]{MR2087303}.}
  \label{thm:rsw_subs_msubs}
  In the notation above, the triples
  \begin{equation*}
  \label{eq:subs_msubs_csp}
  \lp \ch{[0,n - 1]}{k}, C_n, \ch{n}{k}_q \rp \qquad \tx{and} \qquad
      \lp \mch{[0, n - 1]}{k}, C_n, \mch{n}{k}_q \rp
  \end{equation*}
  exhibit the CSP.
\end{thm}

We will also have use of the following principal specializations (see
\cite[Example~I.2.2]{MR1354144} or \cite[Proposition~7.8.3]{MR1676282}):
\begin{align}
  \label{eq:subset_sum}
  \binom{ [0,n - 1]}{k}^{\Sum}(q)
    &= e_k(1, q, q^2, \ldots, q^{n - 1})
     = q^{ \ch{k}{2} } \ch{n}{k}_q, \\
  \label{eq:multiset_sum}
  \mch{ [0,n - 1]}{k}^{\Sum}(q)
    &= h_k(1, q, q^2, \ldots, q^{n - 1})
     = \mch{ n }{k}_q.
\end{align}
Here the $\Sum$ statistic denotes the sum of the elements of a subset or
submultiset of $\bZ$.

Recall that the length function $\ell$ on $S_n$ coincides with the
inversion statistic defined above on words of content
$(1, 1, \ldots, 1)$. More generally, minimal length coset representatives
of parabolic quotients $S_n/S_J$ also have length given by the inversion
statistic on the corresponding words $\W_\alpha$. The following classical
result is due to MacMahon.

\begin{thm}{\cite[Art.~6]{MR1506186}.}
  \label{thm:macmahon}
  For each $\alpha \vDash n$, $\maj$ and $\inv$ are equidistributed on
  $\W_\alpha$ with
\begin{align}
\label{eqn:macmahon}
     \W_\alpha^{\maj}(q) = \binom{n}{\alpha}_q = \W_\alpha^{\inv}(q).
\end{align}
\end{thm}

\section{Modular Periodicity and an Extension Lemma}
\label{modular_periodicity}

We now introduce the concept of \textit{modular periodicity} and use it
to give an extension lemma, Lemma~\ref{lem:extend.csp}, which allows us to
extend CSP's from certain subgroups to larger groups. We
will verify the hypotheses of Lemma~\ref{lem:extend.csp} in the subsequent
sections to deduce Theorem~\ref{thm:alpha_delta}.

\begin{defn}
  \label{defn:period a mod b}
  We say a statistic $ \stat \colon W \to \bZ $ has
  \emph{period $ a $ modulo $ b $ on $ W $} if for all $ i \in \bZ $,
  \[
	 \# \{ w \in W : \stat(w) \equiv_b i \}
     = \# \{ w \in W : \stat(w) \equiv_{b} i + a \}.
  \]
  Similarly, we say a Laurent polynomial $ f(q) \in \bC[q, q^{-1}]$ has
  \emph{period $ a $ modulo $ b $} if
  \[
	q^a f(q) \equiv f(q) \quad (\Mod q^b - 1),
  \]
  or equivalently if $ (q^b - 1) \mid (q^a - 1)f(q) $.
\end{defn}

For example, $ 1 + 5q + q^2 + 5q^3 + q^4 + 5q^5 $ has period 2 modulo $ 6 $. Note that $ \stat $ has period $ a $ modulo $ b $ on $ W $ if and only if
$ W^{\stat}(q) $ has period $ a $ modulo $ b $. The following basic
properties of modular periodicity will be useful throughout the paper.

\begin{lem}
  \label{lem:PeriodFacts}
Let $ f(q) \in \bC[q,q^{-1}] $ and $ a, b, c \in \bZ$.
  \begin{enumerate}[(i)]
    \item If $ f(q) $ has period $ a $ modulo $ c $ and period $ b $ modulo
      $ c $, then $ f(q) $ has period $ u a + v b $ modulo $ c $
      for any $u, v \in \bZ$. In particular, $ f(q) $ has period $ \gcd(a,b) $ modulo $ c $.
    \item If $ f(q) $ has period $ a $ modulo $ b $ and period $ b $ modulo
      $ c $, then $ f(q) $ has period $ a $ modulo $ c $.
    \item If $ f(q) $ has period $ a $ modulo $ c $ and $ b \mid c $, then
      $ f(q) $ has period $ a $ modulo $ b $.
    \item If $ f(q) $ has period $ a $ modulo $ b $, then so does
      $ f(q) h(q) $ for any Laurent polynomial $ h(q) $.
    \item If $ f(q) $ has period $ a $ modulo $ b $ and $ a \mid b $, then
        \[ f(q) \equiv \frac{a}{b} \lp \f{q^b - 1}{q^a - 1} \rp f(q)
           \quad (\Mod q^b - 1). \]
  \end{enumerate}
\end{lem}

\begin{pf} (i), (iii), (iv), and (v) are straightforward. For (ii), suppose
    \[ (q^b - 1) \mid (q^a - 1)f(q), \qquad (q^c - 1) \mid (q^b - 1)f(q). \]
  Write $ q^c - 1 = \prod_{k=1}^c (q - \w_c^k) $. If $q - \w_c^k$ does not divide $f(q)$, then it must divide $q^b - 1$ and hence $q^a - 1$. It follows that
    \[ (q^c - 1) \mid (q^a - 1)f(q). \]
\end{pf}

\begin{lem}\label{lem:extend.csp}
  Suppose $C_n = \langle \sigma_n\rangle$ acts on $W$. Let $g \mid n$ and $C_g \coloneqq \langle \si^{n/g}_n \rangle \subset C_n$. If
  \begin{enumerate}[(i)]
    \item $(W, C_g, f(q))$ exhibits the CSP,
    \item $f(q)$ has period $g$ modulo $n$, and
    \item for all $C_n$-orbits $\cO \subset W$, we have $\frac{n}{|\cO|} \mid
      g$,
  \end{enumerate}
  then $(W, C_n, f(q))$ exhibits the CSP.
\end{lem}

  \begin{pf}
    Let
      \[ F(q) \coloneqq  \sum_{\text{$C_n$-orbits }\cO \subset W } \;
          \frac{q^n - 1}{q^{n/|\cO|} - 1}.
      \]
    By \eqref{CSP2}, $(W, C_n, F(q))$ exhibits the CSP, so
    $(W, C_g, F(q))$ also exhibits the CSP by Remark~\ref{rem:relatedCSPs}.
    Thus, by \eqref{CSP2} and condition (i),
    \begin{align}
    \label{eq:fvsF}
      f(q) = F(q) + p(q)(q^g - 1)
    \end{align}
    for some $ p(q) \in \bC[q] $. Each summand of $ F(q) $ has period
    $g$ modulo $n$ since
      \[ (q^n - 1) \mid (q^g - 1)
         \frac{q^n - 1}{q^{n/|\cO|} - 1}, \]
    by condition (iii). Putting this together with condition (ii),
    $ f(q) $ and $F(q)$ have period $g$ modulo $n$. Using
    Lemma~\ref{lem:PeriodFacts}(v) twice along with \eqref{eq:fvsF} now gives
    \begin{align*}
      f(q)
        &\equiv \frac{g}{n} \frac{q^n - 1}{q^g - 1} f(q) \\
         &= \frac{g}{n} \frac{q^n - 1}{q^g - 1}
            (F(q) + p(q) (q^g - 1) ) \\
         & \equiv \frac{g}{n} \frac{q^n - 1}{q^g - 1} F(q) \\
         & \equiv F(q) \qquad \Mod\text{($q^n - 1$)}.
    \end{align*}
  \end{pf}

\section{Cyclic Descent Type}
\label{sec:CDT}

In this section, we introduce the \textit{cyclic descent type} of a word. We also verify hypothesis (iii) of Lemma~\ref{lem:extend.csp} for
$\W_{\alpha, \delta}$ for a particular $ g $; see Lemma~\ref{lem:freqdivg}.

Let $ w^{(i)} $
denote the subsequence of $ w $
with all letters larger than $ i $ removed. We have a ``filtration''
\[
	\varnothing \preceq w^{(1)} \preceq w^{(2)}
        \preceq \dots \preceq w^{(m - 1)} \preceq w^{(m)} = w,
\]
where $ u \preceq v $ means that $ u $ is a subsequence of $ v $. We think of
this filtration as building up $ w $ by recursively adding all of the copies
of the next largest letter ``where they fit.'' The \textit{cyclic descent
type} of a word $ w $, denoted $ \CDT(w) $, is the sequence which tracks the number of new
cyclic descents at each stage of the filtration. Precisely, we have the following.

\begin{defn}
  \label{def:cdt}
  The \textit{cyclic descent type} (CDT) of a word $w$ is the weak composition of $\cdes(w)$ given by
\begin{align}\label{defn:CDT}
  \begin{split}
	\CDT(w) \coloneqq (\cdes(w^{(1)}), &\cdes(w^{(2)}) - \cdes(w^{(1)}), \dots, \\
                &\cdes(w^{(m)}) - \cdes(w^{(m - 1)}) ).
  \end{split}
\end{align}
\end{defn}

Note that $ \CDT $ is constant on necklaces since rotating $ w $ rotates each
$ w^{(i)} $ and $\cdes$ is constant under rotations. Furthermore, $\cdes(w^{(1)}) = 0$ always, so $\CDT(w)$ always begins with $0$.

\begin{ex}
  \label{ex:cdt}
Suppose $ w = 143124114223 $, so
  \begin{align*}
    & w^{(1)} = 1111              & \cdes(w^{(1)}) = 0, \\
    & w^{(2)} = 112.1122.         & \cdes(w^{(2)}) = 2, \\
    & w^{(3)} = 13.12.11223.      & \cdes(w^{(3)}) = 3, \\
    & w^{(4)} = 14.3.124.114.223. & \cdes(w^{(4)}) = 5.
  \end{align*}
  Hence, $\CDT(143124114223) =
  (0, 2-0, 3-2, 5-3) = (0, 2, 1, 2)$.
\end{ex}

Recall from \eqref{eqn:FixedCDT} that
\begin{align*}
  \W_{\alpha, \delta}
    &\coloneqq \{w \in \W_n : \cont(w) = \alpha, \, \CDT(w) = \delta\}.
\end{align*}
We could define $ \W_{\alpha, \delta} $ more ``symmetrically'' by
replacing $\cont$ with ``cyclic weak ascent type,'' which would
be the point-wise difference of $\cont$ and $\CDT$. However, content is ubiquitous in the literature, so we use it.

\begin{rem}
Despite \eqref{eqn:macmahon}, $ \maj $ and $ \inv $ are not equidistributed even modulo $ n $ on $ \W_{\al,\de} $ in general, so $ (\W_{\al,\de}, C_n, \W_{\al,\de}^{\inv}(q)) $ does not generally exhibit the CSP. For example, $ \W_{(2, 2), (0, 2)} = \{1212, 2121\}$, which has
\[
    \W_{(2, 2), (0, 2)}^{\maj}(q) =  q^2 + q^4, \qquad \W_{(2, 2), (0, 2)}^{\inv}(q) = q^1 + q^3,
\]
which are not even congruent modulo $q^4 - 1$.
\end{rem}

\begin{lem}
\label{lem:freqdivg}
If $\alpha = (\alpha_1, \ldots, \alpha_m)$, $\delta = (\delta_1, \ldots, \delta_m)$, $ N \sub \W_{\al,\de} $ is a necklace, and $ g \coloneqq gcd(\al_1, \dots, \al_m, \de_1, \dots \de_m) $, then $ \f{n}{|N|} \mid g $.

\end{lem}

\begin{pf} Suppose $ N $ is the necklace of $ w $, meaning $ \freq(w) = \f{n}{|N|} $, so we can write $ w = u^{\f{n}{|N|}} $. Hence, using pointwise multiplication,
\[
    \cont(w) = \f{n}{|N|} \dd \cont(u), \qquad \CDT(w) = \f{n}{|N|} \dd \CDT(u).
\]
In particular, $ \f{n}{|N|} $ divides $ \al_1, \dots, \al_m, \de_1, \dots \de_m $, so $ \f{n}{|N|} \mid g $.

\end{pf}

\section{Runs and Falls}
\label{sec:runs}

In this section, we give a method to algorithmically construct
$\W_{\alpha, \delta}$ and use it to prove a product formula for
$\W_{\alpha, \delta}^{\maj}(q)$ modulo $ q^n - 1 $, Theorem~\ref{thm:Maj GF Mod n}. We conclude the
section by using this formula to verify hypothesis (ii) of
Lemma~\ref{lem:extend.csp} for $\W_{\alpha, \delta}$; see Proposition~\ref{Propn:Periodg}.

\subsection{A Tree Decomposition for $\W_{\alpha, \delta}$}

We now describe a way to create words with a fixed content and CDT in terms
of insertions into \textit{runs} and
\textit{falls}. This procedure is organized into a tree, Definition~\ref{def:tree}, whose edges are labeled with sets and multisets.
Lemma~\ref{lem:maj_triple} describes changes in the major index upon traversing an
edge of this tree.

\begin{defn}
  Write $ w = w_1 \cdots w_n \in \W_n $. A \textit{fall} in $w$ is a maximal set of distinct consecutive indices $i, i+1,
  \ldots, j-1, j$ such that $w_i > w_{i+1} > \cdots > w_j$, where we take indices modulo $ n $. A \textit{run}
  in a non-constant word $w$ is a maximal set of distinct consecutive indices $i, i+1,
  \ldots, j$ such that $w_i \leq w_{i+1} \leq \cdots \leq w_j$, where we take indices modulo $ n $. The constant
  word $w = \ell^n$ by convention has no runs and $ n $ falls.
\end{defn}

Note that each letter in $ w $ is part of a unique fall and a unique run, except when $ w = \ell^n $ is constant. It is easy to see that $w$ has $n - \cdes(w)$ falls and $\cdes(w)$ runs, since they are separated by cyclic weak ascents and cyclic descents, respectively. Note that this holds if $ w $ is constant since then $w$ by convention has no runs. Index falls from $0$ from left to right starting at the fall containing the first letter of $w$, and do the same with runs.

\begin{defn}
We write
\begin{equation*}
    F(w) \coloneqq [0, |w| - \cdes(w) - 1]
    \qquad \text{and} \qquad
    R(w) \coloneqq [0,\cdes(w) - 1]
\end{equation*}
for the indices of the falls and runs of $w$, respectively.
\end{defn}

\begin{ex}
  Let $w = 26534611 = 2\udot653\udot4\udot61\udot1\udot = 26.5.346.11 \in \W_8 $,
  where upper dots indicate cyclic weak ascents and lower dots indicate
  cyclic descents. Since $\cdes(w) = 3$, we have
  $F(w) = [0, 4]$ and $R(w) = [0, 2]$. The $5$ falls of $ w $ are
  $ 2, 653, 4, 61, 1 $, with respective indices $0, 1, 2, 3, 4$.
  The $3$ runs of $w$ are $1126, 5, 346$, with
  respective indices $0, 1, 2$.
\end{ex}

\begin{defn} Let $ w $ be a word. Fix a letter $\ell$ and pick a subset $ F $ of the falls $F(w)$. Assume $\ell$ does not appear in any of the falls in $F$. We
  \textit{insert $\ell$ into falls $F$} by successively inserting $\ell$ into
  each fall $w_i > w_{i+1} > \cdots > w_j$ in $ F $ so that $ w_i \cdots \ell \cdots w_j $ is still decreasing.

  Similarly, we may fix a letter $\ell$ and pick a \textit{multisubset} $R$ of $R(w)$ (this time $\ell$ may already appear in a run in $R$). We \textit{insert $\ell$ into runs $R$} by successively inserting $\ell$ into each run $w_i \leq w_{i+1} \leq \cdots \leq w_j$ in $R$ so that $ w_i \cdots \ell \cdots w_j $ is still weakly increasing.
\end{defn}

  When inserting $\ell$ into a run already containing $\ell$, the resulting word is independent of precisely which of the possible positions is used. This is the reason we insert into runs and falls instead of positions.

  Note that there is a slight ambiguity in our description of insertion into falls and runs, since it may be possible to insert either at the beginning or at the end of $w$ while still satisfying the relevant inequalities. Given the choice, we always insert at the beginning of $w$.

\begin{ex}
  Let $w = 2\udot653\udot4\udot61\udot1\udot$\ ~. Insert $7$ into falls of $w$
  with indices $0$ and $3$ to successively obtain
  $\underline{7}2\udot653\udot4\udot61\udot1\udot$ and then $w' \coloneqq
  72\udot653\udot4\udot\underline{7}61\udot1\udot$\ ~. Note that $w'=
  7.26.5.347.6.11$ has two more runs (or cyclic descents) than $w$. Now insert
  $7$ into the runs of $w'$ with multiset of indices $\{ 0,2,3,3 \}$ to successively
  obtain $\underline{7}7.26.5.347.6.11$, $77.26.5\underline{7}.347.6.11$,
  $77.26.57.34\underline{7}7.6.11$, and
  $ w'' \coloneqq 77.26.57.34\underline{7}77.6.11$.
\end{ex}

Let
\begin{align}
    \widetilde{\W}_n & = \{ w \in \W_n : w \tx{ ends in a 1} \}, \\
     \widetilde{\W}_{\al, \de} & = \{ w \in \W_{\al, \de} : w \tx{ ends in a 1} \}.
\end{align}
We restrict to $\widetilde{\W}_n$ and
$\widetilde{\W}_{\alpha, \delta}$ since the major index generating
function is easier to find and extends to $\W_{\alpha, \delta}^{\maj}(q)$
(mod $q^n - 1$).

\begin{defn}
\label{def:FwRw}
Fix $ w \in \widetilde{\W}_n $, a letter $\ell$ not in $ w $, and
\[
     F \subset F(w) = [0, |w| - \cdes(w) - 1] \qquad \text{and} \qquad
     R \underset{\text{mult.}}{\subset} [0, \cdes(w) + |F| - 1]
\]
where $\underset{\text{mult.}}{\subset}$ denotes a multisubset.
Let $w'$ be obtained by inserting $\ell$ into falls $F$ of $w$.
Note that $[0, \cdes(w) + |F|-1] = R(w')$ indexes the runs of $w'$.
Now let $w''$ be obtained by inserting $\ell$ into runs $R$ of $w'$. We
say $w''$ is obtained by \textit{inserting the triple} $(\ell, F, R)$
into $w$. Observe that $\cdes(w'') = \cdes(w') = \cdes(w) + |F|$ and
$w'' \in \wt{W}_{n+|F|+|R|}$.
\end{defn}

We next describe the effect of inserting a single letter on $ \maj $. We restrict to $ \widetilde{\W}_n $ so we preserve a cyclic weak ascent at the end and never add a letter to the end. The fact that the increments in major index from inserting a new letter into all possible positions form a permutation was first observed by Gupta \cite{MR495467}. Lemma~\ref{lem:maj_insertion} tells us exactly the increment in major index based on which run or fall the newly inserted letter fits into.

\begin{lem}
  \label{lem:maj_insertion}
  Suppose $w' \in \widetilde{\W}_{n + 1} $ is obtained by adding a letter $ \ell $ to $ w \in \widetilde{\W}_n $ in any position. Then $w'$ is obtained by inserting $ \ell $ into some run or fall of $ w $, and
\begin{align}
        \maj(w') - \maj(w) = \begin{cases}
    \cdes(w) - r       & \quad \text{if $ \ell $ is inserted into run $ r $ of $ w $}\\
    \cdes(w) + 1 + f  & \quad \text{if $ \ell $ is inserted into fall $ f $ of $ w $.}
  \end{cases}
\end{align}
\end{lem}

\begin{pf} If $\cdes(w') = \cdes(w)$, then $\ell$ is inserted into some run of $w$, and otherwise $\cdes(w') = \cdes(w) + 1$ and $\ell$ is inserted into some fall of $w$. Inserting $\ell$ into run $ r $ of $w$ will increment the position of $ \cdes(w) - r $ descents by $1$ each, so
\begin{align*}
    \maj(w') - \maj(w) = \cdes(w) - r.
\end{align*}
Let $\comaj(w) \coloneqq 1+2+\cdots + (|w| - 1) - \maj(w)$, which is the sum of $i \in [|w| - 1]$ where $w_i \leq w_{i+1}$. Inserting $\ell$ into fall $ f $ of $w$ will increment the position of $ (|w|-1) - \cdes(w) - f $ weak ascents by $1$ each, so
\begin{align*}
    \comaj(w') - \comaj(w) & = (|w|-1) - \cdes(w) - f,
\end{align*}
from which it follows that
\begin{align*}
    \maj(w') - \maj(w) & = \cdes(w) + 1 + f.
\end{align*}

\end{pf}

\begin{lem}
  \label{lem:maj_triple}
  Suppose $w''$ is obtained by inserting the triple $(\ell, F, R)$ into $w \in \widetilde{\W}_n $. Then
    \begin{align}
      \begin{split}
        \maj(w'') - \maj(w) = \binom{|F|+1}{2} &+ (\cdes(w))(|F| + |R|)\\
         &+ |F||R| + \sum_{f \in F} f - \sum_{r \in R} r.
      \end{split}
    \end{align}
\end{lem}

  \begin{pf}
    Let $w'$ be obtained by inserting $\ell$ into falls $F$ of $w$. It
    suffices to show
    \begin{equation}
      \label{eq:maj_diff_1}
      \maj(w') - \maj(w)
        = \binom{|F|+1}{2} + (\cdes(w))|F| + \sum_{f \in F} f
    \end{equation}
    and
    \begin{equation}
      \label{eq:maj_diff_2}
      \maj(w'') - \maj(w')
        = (\cdes(w'))|R| - \sum_{r \in R} r
    \end{equation}
    since $\cdes(w') = \cdes(w'') = \cdes(w) + |F|$. Both \eqref{eq:maj_diff_1} and \eqref{eq:maj_diff_2} follow from iterating Lemma~\ref{lem:maj_insertion} and recalling $ \cdes $ is incremented by 1 each time we insert into a fall.
 \end{pf}

\begin{nota}
\label{not:notation}
For the rest of this section, fix a strong composition $\alpha = (\al_1, \dots, \al_m) $ of $ n \ge 1 $ and $ \delta = (\de_1, \dots, \de_m) \vDash k $ with $\de_1 = 0$. We emphasize that $\alpha$ and $\delta$ have
the same number, $m$, of parts. For $ \ell = 1, \ldots, m $, let
\begin{align}
    n_{\ell} &\coloneqq \alpha_1 + \cdots + \alpha_{\ell}, \\
    k_{\ell} &\coloneqq \delta_1 + \cdots + \delta_{\ell}.
\end{align}
For $ w \in \W_{\al,\de} $, we have the defining conditions $ |w^{(\ell)}| = n_\ell $ and $ \cdes(w^{(\ell)}) = k_\ell $. Furthermore, let
\begin{align*}
    \SSS_\ell &\coloneqq \ch{ [0, n_{\ell - 1} - k_{\ell - 1} - 1]}{\de_\ell},
          \qquad
    \M_\ell \coloneqq \mch{[0, k_\ell - 1]}{\al_\ell - \de_\ell}
  \end{align*}
and
  \[ g \coloneqq \gcd(\alpha_1, \alpha_2, \ldots, \alpha_m, \delta_1, \delta_2, \ldots,
               \delta_m). \]
\end{nota}

If $w \in \widetilde{\W}_{\alpha, \delta}$, then the set $\SSS_\ell$ consists
of all subsets of the falls $F(w^{(\ell-1)})$ which, when $\ell$ is
inserted into those falls of $w^{(\ell-1)}$, result in a word $w'$ with
$k_\ell$ cyclic descents. The multiset $\M_\ell$ similarly consists of all
choices of runs $R(w')$ which, when $\ell$ is inserted into those runs,
result in a word with length $n_\ell$.

\begin{rem}
  We restrict to strong compositions $\alpha$ for notational
  simplicity, though the results in this section may easily be
  generalized to arbitrary weak compositions by ``flattening'' weak compositions to strong ones by removing zeros.
\end{rem}

\begin{defn}
  \label{def:tree}
  Construct a rooted, vertex-labeled and edge-labeled tree
  $T_{\alpha,\delta}$ recursively as follows. Begin with a tree $ T^{(1)} $ containing only a root labeled by the word $1^{\al_1}$. For $\ell=2, \ldots, m$,
  to obtain $T^{(\ell)}$, do the following. For each leaf $w$ of
  $ T^{(\ell - 1)} $ and for each triple $(\ell, F, R)$ with
   \[ F \in \SSS_\ell \qquad \text{and} \qquad
         R \in \M_\ell, \]
add an edge labeled by $(F, R)$ to $ T^{(\ell - 1)} $ from $ w $ to $w''$ where $w''$ is obtained by inserting $(\ell, F, R)$ into $w$. Define $ T_{\al, \de} \coloneqq T^{(m)} $.
\end{defn}

\begin{ex}
  \label{ex:tree}
  Let $ \alpha = (3,1,1) $ and $ \delta = (0,1,0) $.
  Figure~\ref{fig:tree1} is the tree $T_{\alpha,\delta}$.

  \begin{figure}[ht]
      \centering
    \begin{tikzpicture}
      \node (2)  at (0, -4)  {$111$};
      \node (3a) at (-3.5, -6) {$2111$};
      \node (3b) at (0, -6)    {$1211$};
      \node (3c) at (3.5, -6)  {$1121$};
      \node (4a) at (-3.5, -8) {$23111$};
      \node (4b) at (0, -8)    {$12311$};
      \node (4c) at (3.5, -8)  {$11231$};

      \draw [->] (2) -- (3a) node [pos = 0.3, left]
        {$(\{0 \}, \varnothing)$};
      \draw [->] (2) -- (3b) node [pos = 0.6, right]
        {$(\{1 \}, \varnothing)$};
      \draw [->] (2) -- (3c) node [pos = 0.3, right]
        {$( \{2\}, \varnothing )$};
      \draw [->] (3a) -- (4a) node [midway, left]
        {$(\varnothing, \{ 0 \})$};
      \draw [->] (3b) -- (4b) node [midway, right]
        {$(\varnothing, \{ 0 \})$};
      \draw [->] (3c) -- (4c) node [midway, right]
        {$(\varnothing, \{ 0 \})$};
    \end{tikzpicture}
    \caption{Full tree $T_{\alpha, \delta}$ with
      $\alpha = (3,1,1)$, $\delta = (0,1,0)$.}
      \label{fig:tree1}
  \end{figure}
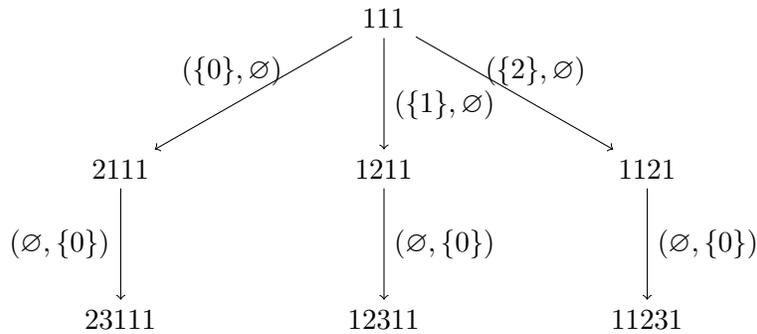

  Let $\alpha = (4, 2, 3)$ and $\delta = (0, 2, 1)$. Figure~\ref{fig:tree2} is the
  subgraph of $T_{\alpha,\delta}$ consisting of paths from the root to leaves
  that are rotations of $112113323$. For this full $T_{\alpha, \delta}$, the root has $\binom{4}{2} = 6$
  children since $1111$ has $4$ falls. Each child of the root itself has
  $\binom{4}{1}\mch{3}{2} = 24$ children. Hence, $T_{\alpha,\delta}$ has
  $144$ leaves. Notice that the cyclic rotations of $311211332$ appearing
  as leaves in Figure~\ref{fig:tree2} are precisely those ending in $1$.
  It will shortly become apparent that in this example,
  $\#\W_{\alpha, \delta} = \frac{9}{4} \cdot 144 = 324$.

  \begin{figure}[ht]
      \centering
    \begin{tikzpicture}
      \node (2)  at (0, -4.5)  {$1111$};
      \node (3a) at (-2.5, -6) {$211211$};
      \node (3b) at (2.5, -6)  {$121121$};
      \node (4a) at (-4, -8) {$332311211$};
      \node (4b) at (-1.5, -8) {$211332311$};
      \node (4c) at (1.5, -8)  {$133231121$};
      \node (4d) at (4, -8)  {$121133231$};

      \draw [->] (2) -- (3a) node [midway, left]
        {$(\{0, 2\}, \varnothing)$};
      \draw [->] (2) -- (3b) node [midway, right]
        {$(\{1, 3\}, \varnothing)$};
      \draw [->] (3a) -- (4a) node [midway, left]
        {$( \{0\}, \{ 0, 1 \})$};
      \draw [->] (3a) -- (4b) node [auto = left, pos = 0.35 ]
        {$(\{2\}, \{ 1, 2 \})$};
      \draw [->] (3b) -- (4c) node [auto = right, pos = 0.9]
        {$(\{1\}, \{ 0, 1 \})$};
      \draw [->] (3b) -- (4d) node [midway, right]
        {$(\{3\}, \{ 1, 2\})$};
    \end{tikzpicture}
    \caption{Subgraph of tree $T_{\alpha, \delta}$ with
      $\alpha = (4, 2, 3)$, $\delta = (0, 2, 1)$.}
      \label{fig:tree2}
  \end{figure}
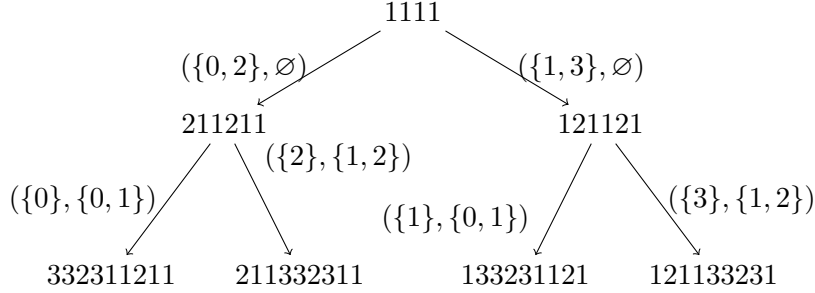

\end{ex}

\begin{lem}
  \label{lem:T_alpha_delta}
  The vertices of $T_{\alpha, \delta}$ which are $ \ell < m $ edges away from
  the root are precisely the elements of $\{w^{(\ell + 1)} : w \in
  \widetilde{\W}_{\alpha, \delta} \}$, each occurring once. In particular, the
  leaves of $T_{\alpha, \delta}$ are precisely the elements of
  $\widetilde{\W}_{\alpha, \delta}$, each occurring once.
\end{lem}

  \begin{pf}
   By definition of $ \SSS_\ell $ and $ \M_\ell $, any leaf of $ T^{(\ell)} $
   has content $ (\al_1, \dots, \al_\ell) $, cyclic descent type
   $ (\de_1, \dots, \de_\ell) $, and ends in a 1, so is in
   $\{w^{(\ell)} : w \in
  \wt{\W}_{\alpha, \delta} \}$. Conversely, given any
   $w \in \widetilde{\W}_{\alpha, \delta}$, the word $w^{(\ell)}$ is obtained
   by inserting a unique triple $(\ell, F, R)$ into $w^{(\ell-1)}$ by repeated applications
   of Lemma~\ref{lem:maj_insertion}.
  \end{pf}

\begin{defn}
  \label{def:Phi}
  By Lemma~\ref{lem:T_alpha_delta}, the tree $T_{\alpha, \delta}$ encodes a
  bijection
    \[ \Phi \colon \widetilde{\W}_{\alpha, \delta}
       \too{\sim} \prod_{\ell=2}^m \SSS_\ell \times \M_\ell \]
  given by reading the edge labels from the root to $w$. We suppress the dependence of $\Phi$ on $\alpha$ and $\delta$ from the notation since they can be computed from the input $w$.
\end{defn}

\begin{lem}
\label{lem:tildefraction}
  For any $w \in \W_{\alpha, \delta}$,
  \begin{equation}
    \label{eq:twiddle_fraction_necklace}
    \#[w] =
    \frac{n}{\alpha_1} \cdot
    \# \lp [w] \cap \widetilde{\W}_{\alpha, \delta} \rp.
  \end{equation}
  Consequently,
  \begin{equation}
    \label{eq:twiddle_fraction}
    \#\W_{\alpha, \delta} = \frac{n}{\al_1} \cdot \#\widetilde{\W}_{\alpha, \delta}.
  \end{equation}
\end{lem}

  \begin{pf}
    Each $w \in \W_{\alpha, \delta}$ has $\period(w) = n/\freq(w)$
    distinct cyclic rotations, of which $\alpha_1/\freq(w)$ end in $1$.
  \end{pf}

\begin{prop}
  \label{prop:W_alpha_delta_count}
  Using Notation~\ref{not:notation}, we have
  \begin{equation}
    \label{eq:W_alpha_delta}
    \# \W_{\alpha, \delta}
       = \frac{n}{\alpha_1}
         \prod_{\ell=2}^{m} \binom{n_{\ell-1} - k_{\ell-1}}{\delta_\ell}
           \mch{k_\ell}{\alpha_\ell - \delta_\ell}.
  \end{equation}
  In particular, $\W_{\alpha, \delta} \neq \varnothing$ if and only if
  \begin{equation}
    \label{eq:polytope}
    \begin{aligned}
      0 &\leq \delta_\ell \leq \alpha_\ell
        &\text{for all }1 \leq \ell \leq m, \text{and} \\
      \delta_1 + \cdots + \delta_{\ell+1}
        &\leq \alpha_1 + \cdots + \alpha_\ell
        &\text{for all }1 \leq \ell < m.
    \end{aligned}
  \end{equation}
\end{prop}

  \begin{pf}
    The product in \eqref{eq:W_alpha_delta} is
    $\# \prod_{\ell=2}^m \SSS_\ell \times \M_\ell$, which is
    $\#\widetilde{\W}_{\alpha, \delta}$ by the bijection $\Phi$. Now
    \eqref{eq:W_alpha_delta} follows from \eqref{eq:twiddle_fraction}, and \eqref{eq:polytope} follows from \eqref{eq:W_alpha_delta}.
  \end{pf}

\subsection{Major Index Generating Functions}

We next use the bijection $\Phi$ and Lemma~\ref{lem:maj_triple} to give a product formula for $\widetilde{\W}_{\alpha, \delta}^{\maj}(q) $,
Theorem~\ref{thm:tW_alpha_delta_maj}. We then use modular periodicity to obtain an analogous expression for $ \W_{\al,\de}^{\maj}(q) $ modulo $ q^n - 1 $, Theorem~\ref{thm:Maj GF Mod n}.

\begin{thm} Using Notation~\ref{not:notation}, we have
  \label{thm:tW_alpha_delta_maj}
    \begin{align}
    \label{eq:tW_alpha_delta_maj}
    \widetilde{\W}_{\alpha, \delta}^{\maj}(q)
         & = \prod_{\ell=2}^{m} q^{k_\ell \al_\ell}
         \binom{n_{\ell - 1} - k_{\ell - 1}}{\delta_{\ell}}_q
         \mch{k_{\ell}}{\alpha_{\ell} - \delta_{\ell}}_{q^{-1}} \\
    \label{eq:tW_alpha_delta_maj2}
         & = q^{\eta(\al,\de)} \prod_{\ell=2}^{m}
         \binom{n_{\ell - 1} - k_{\ell - 1}}{\delta_{\ell}}_q
         \mch{k_{\ell}}{\alpha_{\ell} - \delta_{\ell}}_{q}
    \end{align}
  where
    \[ \eta(\al,\de) \coloneqq n - \alpha_1 + \binom{k}{2}
               + \sum_{\ell=2}^m \binom{\delta_\ell}{2}. \]
\end{thm}

  \begin{pf}
Recall $ \Sum(A) $ denotes the sum of the elements of a set or multiset $ A $. Combining $\Phi$ with Lemma~\ref{lem:maj_triple} shows that
    \begin{align}
      \label{eq:W_tilde_lFR}
      \widetilde{\W}_{\alpha, \delta}^{\maj}(q)
        &= \prod_{\ell=2}^{m}
           \sum_{\substack{F \in \SSS_\ell \\ R \in \M_\ell}}
             q^{\epsilon(\ell, F, R)},
    \end{align}
    where
     \[ \epsilon(\ell, F, R) \coloneqq \binom{\delta_{\ell}+1}{2} + k_{\ell - 1} \alpha_{\ell}
                + \delta_{\ell} (\alpha_{\ell} - \delta_{\ell})
                + \Sum(F) - \Sum(R). \]
   Noting that
      \[ \binom{\delta_\ell+1}{2} + k_{\ell-1}\alpha_\ell
         + \delta_\ell(\alpha_\ell - \delta_\ell)
         = k_\ell \alpha_\ell - \binom{\delta_\ell}{2},
      \]
    simplifying \eqref{eq:W_tilde_lFR} gives
    \begin{align}
    \label{eqn:majGFsimp1}
      \widetilde{\W}_{\alpha, \delta}^{\maj}(q)
        &= \prod_{\ell=2}^{m} q^{k_{\ell}\alpha_{\ell} - \binom{\delta_{\ell}}{2}}
             \SSS_\ell^{\Sum}(q)
             \M_\ell^{\Sum}
             (q^{-1}).
    \end{align}
    Equation \eqref{eq:tW_alpha_delta_maj} now follows from
    \eqref{eq:subset_sum}, \eqref{eq:multiset_sum}, and
    the definition of $\SSS_\ell$ and $\M_\ell$. As for
    \eqref{eq:tW_alpha_delta_maj2}, consider the reversal bijection $ r \colon \M_\ell \to \M_{\ell} $ induced by
      \[ x \mapsto k_\ell - 1 - x \]
    on $[0, k_\ell-1]$.
    This bijection satisfies $ \Sum(r(A)) = (k_\ell - 1)(\al_\ell - \de_\ell) - \Sum(A) $, so
\begin{align}
\label{eqn:sumreflected}
    \M_\ell^{\Sum}(q^{-1})
      = q^{-(k_\ell - 1)(\al_\ell - \de_\ell)} \M_\ell^{\Sum}(q).
\end{align}
Plugging \eqref{eqn:sumreflected} into \eqref{eqn:majGFsimp1} and noting
that
\begin{align*}
  \sum_{\ell=2}^m \left(k_\ell \alpha_\ell - \binom{\delta_\ell}{2}
       - (k_\ell - 1)(\alpha_\ell - \delta_\ell)\right)
    &= \sum_{\ell=2}^m \left(\alpha_\ell -\frac{\delta_\ell}{2} - \frac{\delta_\ell^2}{2} + k_\ell \delta_\ell\right) \\
    &= n-\alpha_1 - \frac{k}{2} + \sum_{\ell=2}^m \left(-\frac{\delta_\ell^2}{2} + \sum_{j=2}^\ell \delta_j \delta_\ell\right) \\
    &= n-\alpha_1 - \frac{k}{2} + \frac{1}{2}\sum_{\ell=2}^m \sum_{j=2}^m \delta_j \delta_\ell \\
    &= n-\alpha_1 - \frac{k}{2} + \frac{k^2}{2}
\end{align*}
gives
    \begin{align*}
      \widetilde{\W}_{\alpha, \delta}^{\maj}(q)
        &= q^{n - \al_1 + \binom{k}{2}} \prod_{\ell=2}^{m}
             \SSS_\ell^{\Sum}(q) \M_\ell^{\Sum}(q).
    \end{align*}

  \noindent Using \eqref{eq:subset_sum} and \eqref{eq:multiset_sum} now
  yields \eqref{eq:tW_alpha_delta_maj2}.
  \end{pf}

\begin{lem}
  \label{lem:maj_period}
  Let $\alpha \vDash n$, $\delta \vDash k$.
  The statistic $\maj$ has period $k$ modulo $n$ on $\W_{\alpha, \delta}$.
  Moreover, $\maj$ is constant modulo $d \coloneqq \gcd(n, k)$ on necklaces in
  $\W_{\alpha, \delta}$, and
  \begin{equation}
    \label{eq:W_twiddle_W_maj}
    W^{\maj}_{\alpha, \delta}(q) \equiv \frac{n}{\alpha_1}
       \widetilde{\W}_{\alpha, \delta}^{\maj}(q)
       \qquad \text{(mod $q^d-1$)}.
  \end{equation}
\end{lem}

  \begin{pf}
    Since cyclically rotating $w \in \W_{\alpha, \delta}$ increments each
    cyclic descent by $1$ modulo $n$, we have
    \begin{align}
    \label{eq:majnecklace}
      \maj(\sigma_n \cdot w) \equiv_n \maj(w) + k.
    \end{align}
    In particular, $\maj$ has period $k$ modulo $n$ on necklaces in $ \W_{\al,\de} $. Furthermore, $\maj$ is
    constant on necklaces in $\W_{\al, \de}$ modulo $d$. By
    \eqref{eq:twiddle_fraction_necklace}, each necklace $ [w] \in \W_{\al,\de} $ has the same fraction, $ \f{\al_1}{n} $, of its elements in $ \widetilde{\W}_{\al,\de} $, so \eqref{eq:W_twiddle_W_maj} follows.
  \end{pf}

\begin{thm}
\label{thm:Maj GF Mod n}
  Using Notation~\ref{not:notation}, let $ d \coloneqq \gcd(n, k) $. Then, modulo
  $ q^n - 1 $,
  \begin{align}
  \begin{split}
      \W_{\al,\de}^{\maj}(q) & \equiv \f{d}{\al_1} \lp \f{q^n - 1}{q^d - 1} \rp
         \prod_{\ell=2}^m q^{k_\ell \al_\ell}
           \ch{n_{\ell - 1} - k_{\ell - 1}}{\de_\ell}_q
           \mch{k_\ell}{\al_\ell - \de_\ell}_{q^{-1}} \\
           & \equiv \f{d}{\al_1} \lp \f{q^n - 1}{q^d - 1} \rp
         q^{\ch{k}{2} + \sum_{\ell = 2}^m \ch{\de_\ell}{2} - \al_1 }
         \prod_{\ell = 2}^m \ch{n_{\ell - 1} - k_{\ell - 1}}{\de_\ell}_q
           \mch{k_\ell}{\al_\ell - \de_\ell}_q.
  \end{split}
  \end{align}
\end{thm}

\begin{pf}
  By Lemma~\ref{lem:maj_period}, $\maj$ has period $k$ modulo $n$ on
  $\W_{\alpha, \delta}$. Hence by Lemma~\ref{lem:PeriodFacts}(i), $\maj$ has
  period $d$ modulo $n$ on $\W_{\alpha, \delta}$. Using
  Lemma~\ref{lem:PeriodFacts}(v) and \eqref{eq:W_twiddle_W_maj} gives
  \begin{align*}
    \W_{\alpha, \delta}^{\maj}(q)
      &\equiv \frac{d}{n} \lp \f{q^n - 1}{q^d - 1} \rp
        \W_{\alpha, \delta}^{\maj}(q) \\
      &\equiv \frac{d}{n} \lp \f{q^n - 1}{q^d - 1} \rp \lp \frac{n}{\alpha_1} \widetilde{\W}_{\alpha, \delta}^{\maj}(q) + p(q)(q^d - 1) \rp \\
      & \equiv \f{d}{\al_1} \lp \f{q^n - 1}{q^d - 1} \rp \widetilde{\W}_{\alpha, \delta}^{\maj}
        \qquad (\Mod q^n - 1),
  \end{align*}
where $p(q) \in \bC[q]$. Theorem~\ref{thm:Maj GF Mod n} now follows from
  Theorem~\ref{thm:tW_alpha_delta_maj}.
\end{pf}

\begin{cor}
  \label{cor:vandermonde}
  Using Notation~\ref{not:notation}, let $ d \coloneqq \gcd(n, k) $. Then, modulo
  $q^n - 1$,
  \begin{align}
  \begin{split}
      \W_\alpha^{\maj}(q)
         &= \binom{n}{\alpha}_q
         \equiv \sum_{\delta} \f{d}{\al_1} \lp \f{q^n - 1}{q^d - 1} \rp
	     \prod_{\ell=2}^m q^{k_\ell \al_\ell}
           \ch{n_{\ell - 1} - k_{\ell - 1}}{\de_\ell}_q
	       \mch{k_\ell}{\al_\ell - \de_\ell}_{q^{-1}}
  \end{split}
  \end{align}
  where the sum is over weak compositions $\delta$ of $k$ satisfying
  \eqref{eq:polytope}. In particular,
  \begin{align}
  \label{eq:vandermonde}
  \begin{split}
      \#\W_\alpha
         &= \binom{n}{\alpha}
         = \sum_{\delta} \f{n}{\al_1}
	     \prod_{\ell=2}^m
           \ch{n_{\ell - 1} - k_{\ell - 1}}{\de_\ell}
	       \mch{k_\ell}{\al_\ell - \de_\ell}.
  \end{split}
  \end{align}
\end{cor}

Note that the two-letter case of \eqref{eq:vandermonde} is a special
case of the classical Vandermonde convolution identity
\cite[Ex.~1.1.17]{MR2868112}.

\subsection{Verifying Hypothesis (ii) of Lemma~\ref{lem:extend.csp} for \texorpdfstring{$\W_{\alpha, \delta}$}{Words of Fixed Content and CDT}}

\begin{prop}\label{Propn:Periodg}
  Using Notation~\ref{not:notation}, $ \W_{\al,\de}^{\maj}(q)$ has period $ g $ modulo $ n $.
\end{prop}

\begin{pf}
Let $ d = \gcd(n,k) $. By Theorem~\ref{thm:Maj GF Mod n},
  \begin{align*}
    \W_{\al,\de}^{\maj}(q) \equiv \f{d}{\al_1} \lp \f{q^n - 1}{q^d - 1} \rp
        \prod_{\ell=2}^m q^{k_\ell \al_\ell}
        \ch{n_{\ell - 1} - k_{\ell - 1}}{\de_\ell}_q
        \mch{k_\ell}{\al_\ell - \de_\ell}_{q^{-1}}
  \end{align*}
  modulo $ q^n - 1 $. The action of rotation on elements of
  $ \SSS_\ell = \binom{[0, n_{\ell-1} - k_{\ell-1}]}{\delta_\ell} $
  increases their sum by $ \de_\ell $ modulo $ n_{\ell-1} - k_{\ell-1} $.
  Thus by
  \eqref{eq:subset_sum}, $\ch{n_{\ell - 1} - k_{\ell - 1}}{\de_\ell}_{q}$
  has period $\de_\ell$ modulo $n_{\ell - 1} - k_{\ell - 1}$.
  Similarly by \eqref{eq:multiset_sum},
  $\mch{k_{\ell}}{\al_{\ell} - \de_{\ell}}_{q^{-1}}$ has period
  $ \al_{\ell} - \de_{\ell} $ modulo $ k_{\ell}$.
  For $\ell=2, \ldots, m$, by Lemma~\ref{lem:PeriodFacts}(iv) we then have
  \begin{align}
  \label{eq:perfact1}
    & \tx{ $\W_{\al,\de}^{\maj}(q)$ has period $\de_\ell$
      modulo $n_{\ell - 1} - k_{\ell - 1}$, and } \\
    \label{eq:perfact2}
    & \tx{ $\W_{\al,\de}^{\maj}(q)$ has period $\al_{\ell} - \de_{\ell}$
      modulo $k_{\ell}$.}
  \end{align}

We show $\W_{\al,\de}^{\maj}(q)$ has period $\alpha_\ell$ and
  $\delta_\ell$ modulo $n$ by downward induction on $\ell$, for
  $m \geq \ell \geq 2$. Note that the base case $ \ell = m $ is accounted
  for by our argument as well.

  Suppose $\W_{\al,\de}^{\maj}(q)$ has period $\alpha_j$ and
  $\delta_j$ modulo $n$ for all $j > \ell$. By Lemma~\ref{lem:maj_period},
  $\W_{\alpha, \delta}^{\maj}(q)$ has period $k$ modulo $n$. By
  Lemma~\ref{lem:PeriodFacts}(i),
  $\W_{\al,\de}^{\maj}(q)$ thus has period
    \[ k_\ell = k - (\delta_m + \cdots + \delta_{\ell+1}) \]
  modulo $n$. Since $\W_{\al,\de}^{\maj}(q)$ has period $\al_{\ell} - \de_{\ell}$ modulo $k_{\ell}$, $\W_{\al,\de}^{\maj}(q)$ has period
  $\alpha_\ell-\delta_\ell$ modulo $n$ by Lemma~\ref{lem:PeriodFacts}(ii).

  As noted, $\W_{\al,\de}^{\maj}(q)$ has period $\delta_\ell$ modulo
  $n_{\ell-1} - k_{\ell-1}$. By Lemma~\ref{lem:PeriodFacts}(i),
  $\W_{\al,\de}^{\maj}(q)$ also has period
    \[ n_{\ell-1} - k_{\ell-1}
       = n - (\alpha_m + \cdots + \alpha_{\ell+1})
       - k + (\delta_m + \cdots + \delta_{\ell+1})
       - (\alpha_\ell - \delta_\ell)
    \]
  modulo $n$. Hence, as $\W_{\al,\de}^{\maj}(q)$ has period $\de_\ell$
      modulo $n_{\ell - 1} - k_{\ell - 1}$, $\W_{\al,\de}^{\maj}(q)$
  has period $\delta_\ell$ modulo $n$ by Lemma~\ref{lem:PeriodFacts}(ii). By another application of
  Lemma~\ref{lem:PeriodFacts}(i), $\W_{\al,\de}^{\maj}(q)$ has period
  $\alpha_\ell$ modulo $n$ as well, completing the induction.

  Indeed, $\W_{\al,\de}^{\maj}(q)$ has period $\delta_1 = 0$
  modulo $n$ trivially, and $\W_{\al,\de}^{\maj}(q)$ has period
  $\alpha_1 = n - (\alpha_m + \cdots + \alpha_2)$ modulo $n$ by
  Lemma~\ref{lem:PeriodFacts}(i). Putting everything together, $ \W_{\al,\de}^{\maj}(q) $ has periods $ \al_1, \dots, \al_m, \de_1, \dots, \de_m $ modulo $ n $, so by one more application of Lemma~\ref{lem:PeriodFacts}(i), $ \W_{\al,\de}^{\maj}(q) $ has period $ g $ modulo $ n $.
\end{pf}

\section{Refining the CSP to fixed content and Cyclic Descent Type}
\label{sec:CSP}

In this section, we verify the final hypothesis (i) of Lemma~\ref{lem:extend.csp}
for $\W_{\alpha, \delta}$ and deduce Theorem~\ref{thm:alpha_delta}. Throughout
this section we continue to follow Notation~\ref{not:notation}. We recall in
particular that
\begin{align*}
    \SSS_\ell &\coloneqq \ch{ [0, n_{\ell - 1} - k_{\ell - 1} - 1]}{\de_\ell},
          \qquad
    \M_\ell \coloneqq \mch{[0, k_\ell - 1]}{\al_\ell - \de_\ell}
\end{align*}
and
  \[ g \coloneqq \gcd(\alpha_1, \ldots, \alpha_m, \delta_1, \ldots, \delta_m).
  \]

\subsection{A Fixed Point Lemma}
\label{subsec:fixed_point}

To prove our main result, Theorem~\ref{thm:alpha_delta}, one approach would be
to find a $C_n$-equivariant isomorphism between a known CSP triple and
$(\W_{\alpha, \delta}, C_n, \W_{\alpha, \delta}^{\maj}(q))$. Such a triple
is hinted at by \eqref{eq:tW_alpha_delta_maj} and the bijection $\Phi$
using products of CSP's coming from Theorem~\ref{thm:rsw_subs_msubs}, though the
approach encounters immediate
difficulties. For instance, $\widetilde{\W}_{\alpha, \delta}$ is not
generally closed under the $C_n$-action. In this section, we instead give
a fixed point lemma, Lemma~\ref{lem:fix_iff}, which is intuitively a weakened
version of the equivariant isomorphism approach.

\begin{defn}
  \label{def:restricted_actions}
  We define $C_g$-actions on $\SSS_\ell$, $\M_\ell$, and $\W_{\alpha, \delta}$
  as follows. Since $g \mid n_{\ell-1} - k_{\ell-1}$, $g \mid k_\ell$, and
  $g \mid n$, $C_g$ acts on each of $\SSS_\ell$, $\M_\ell$, and $\W_{\alpha,
  \delta}$ by restricting the actions of
  $C_{n_{\ell-1}-k_{\ell-1}}$, $C_{k_\ell}$, and $C_n$ to their unique
  subgroups of size $ g $. For instance, the action of $C_g$ on
  $\W_{\alpha, \delta}$ is generated by rotation by $n/g$.

  We additionally define $C_g$-actions on $\SSS_\ell \times \M_\ell$ and
  $\prod_{\ell=2}^m \SSS_\ell \times \M_\ell$ by letting $C_g$ act diagonally.
  We emphasize that despite having $C_g$-actions on $\W_{\alpha, \delta}$
  and $\prod_{\ell=2}^m \SSS_\ell \times \M_\ell$, the bijection
  $\Phi \colon \widetilde{\W}_{\alpha, \delta} \too{\sim}
  \prod_{\ell=2}^m \SSS_\ell \times \M_\ell$ is not in general equivariant
  since $ \widetilde{\W}_{\alpha, \delta} $ is not closed under the
  $ C_g $ action on $ \W_{\al,\de} $.
\end{defn}

\begin{defn}
\label{defn:multword}
  Given a multisubset of some set $[0, a]$, we may encode
  it as a multiplicity word $w_0 w_1 \ldots w_a$ where
  $w_i$ is the multiplicity of $i$. In particular, we may consider the
  bijection $\Phi \colon \widetilde{\W}_{\alpha, \delta} \too{\sim}
  \prod_{\ell=2}^m \SSS_\ell \times \M_\ell$ as mapping words to sequences
  of pairs of certain words.
\end{defn}

\begin{ex}
  Consider the leaf $w = 211332311$ in Figure~\ref{fig:tree2} from Example~\ref{ex:tree}.
  Reading edge
  labels gives $\Phi(w) = ((\{0, 2\}, \varnothing), (\{2\}, \{1,2\}))$.
  Recalling that $\SSS_2$ consists of subsets of $[0, 4-1]$, $\M_2$
  consists of multisubsets of $\varnothing$, $\SSS_3$ consists of subsets of
  $[0, 4-1]$, and $\M_3$ consists of multisubsets of $[0, 3-1]$, the
  corresponding sequence of words is $((1010, \epsilon), (0010, 011))$,
  where $\epsilon$ denotes the empty word. Table~\ref{tab:words} summarizes
  several similar translations.

  \begin{table}[ht]
    \begin{tabular}[c]{c|c|c}
      $w$ & $\Phi(w)$ & sequence of pairs of words \\
      \hline
      $211332311$ & $((\{0, 2\}, \varnothing), (\{2\}, \{ 1,2 \}))$
                  & $((1010, \epsilon), (0010, 011))$ \\
      \hline
      $121133231$ & $((\{1, 3\}, \varnothing), (\{3\}, \{ 1,2 \}))$
                  & $((0101, \epsilon), (0001, 011))$ \\
      \hline
      $(211332311)^2$
                  & $((\{0, 2, 4, 6\}, \varnothing)$,
                  & $((1010^2, \epsilon), (0010^2, 011^2))$ \\
     \ & $(\{2, 6\}, \{ 1,2, 4, 5 \} ))$ & \ \\
     \hline
     2221123311 & $((\{0,2\}, \{0,0\}), (\varnothing, \{1,1\}))$
                  & $((1010,20),(000000,02))$
    \end{tabular}
    \vspace{1em}
    \caption{Examples of words, corresponding sequences of edge labels
      in $T_{\alpha, \delta}$, and corresponding sequences of words.
      Note that the second word is a cyclic rotation of the first.}
    \label{tab:words}
  \end{table}
\end{ex}

\begin{lem}
  \label{lem:phis_monoid}
  Suppose $w = u^k$ for some word $u$. If
    \[ \Phi(u) = ((x_2, y_2), \ldots, (x_m, y_m)) \]
  encoded as multiplicity words as in Definition~\ref{defn:multword}, then
    \[ \Phi(w) = ((x_2^k, y_2^k), \ldots, (x_m^k, y_m^k)). \]
\end{lem}

  \begin{pf}
    The insertion triples needed to build $ w $ are the sequences of $ k $ shifted copies of the insertion triples needed to build $ u $.
  \end{pf}

\begin{lem}
  \label{lem:fix_iff}
  An element $\tau \in C_g$ fixes $w \in \widetilde{\W}_{\alpha, \delta}$
  if and only if $\tau$ fixes $\Phi(w)$.
\end{lem}

  \begin{pf}
    For $ \tau \in C_n $, let $ o(\tau) $ denote the order of $ \tau $.
    It is easy to see that $\tau \in C_n$ fixes $w \in \W_n$ if and
    only if there is some word $u$ such that $w = u^{o(\tau)}$.

    Suppose $\tau \in C_g$ fixes $w$, so that
    $w = u^{o(\tau)}$. By Lemma~\ref{lem:phis_monoid},
      \[ \Phi(w) = ((x_2^{o(\tau)}, y_2^{o(\tau)}), \ldots,
         (x_m^{o(\tau)}, y_m^{o(\tau)})). \]
Each of the words $x_i^{o(\tau)}$ and $y_i^{o(\tau)}$
    is fixed by $\tau$, so $\Phi(w)$ is fixed by $\tau$. The reverse
    implication follows analogously using the fact that $\Phi$ is a
    bijection.
  \end{pf}

\subsection{Verifying Hypothesis (i) of Lemma~\ref{lem:extend.csp} for \texorpdfstring{$\W_{\alpha, \delta}$}{Words of Fixed Content and CDT}}

\begin{thm}\label{thm:CSPmodg}
  Using Notation~\ref{not:notation}, $ (\W_{\al, \de}, C_g, \W_{\al, \de}^{\maj}(q)) $ exhibits the CSP.
\end{thm}

  \begin{pf}
    We use the notation and actions in Definition~\ref{def:restricted_actions}.
    Recall that
      \[ \SSS_\ell \coloneqq \binom{[0, n_{\ell-1} - k_{\ell-1} - 1]}{\delta_\ell},
         \qquad
         \M_\ell \coloneqq \mch{[0, k_\ell-1]}{\alpha_\ell - \delta_\ell}. \]
    From Theorem~\ref{thm:rsw_subs_msubs}, for each $ 2 \le \ell \le m $,
      \[ \lp \SSS_\ell, C_g, \ch{n_{\ell} - k_{\ell}}{\de_\ell}_q \rp
         \qquad\text{and}\qquad
         \lp \M_\ell, C_g, \mch{k_\ell}{\al_\ell - \de_\ell}_{q^{-1}} \rp
      \]
    exhibit the CSP. Taking products,
    \begin{align*}
      \lp \prod_{\ell = 2}^m \SSS_\ell \times \M_\ell, C_g,
        \prod_{\ell = 2}^{m} \ch{n_{\ell} - k_{\ell}}{\de_\ell}_q
        \mch{k_\ell}{\al_\ell - \de_\ell}_{q^{-1}} \rp
    \end{align*}
    exhibits the CSP. Comparing this to Theorem~\ref{thm:tW_alpha_delta_maj},
    we have
    \begin{align*}
      \widetilde{\W}^{\maj}_{\alpha, \delta} \equiv \prod_{\ell = 2}^{m}
        \ch{n_{\ell} - k_{\ell}}{\de_\ell}_q
        \mch{k_\ell}{\al_\ell - \de_\ell}_{q^{-1}}
    \end{align*}
    modulo $ q^{g} - 1 $, as $\sum_{\ell=2}^m k_\ell \alpha_\ell
    \equiv_g 0 $ because $ g \mid \al_{\ell} $ for all $ \ell $. Thus,
    \begin{align}
    \label{eqn:CSPSMCW}
	   \lp \prod_{\ell = 2}^m \SSS_\ell \times \M_\ell, C_g,
         \widetilde{\W}_{\al, \de}^{\maj}(q) \rp
    \end{align}
    exhibits the CSP.

By Lemma~\ref{lem:tildefraction}, for any
    $w \in \W_{\alpha, \delta}$,
    \begin{equation*}
      \#[w] = \frac{n}{\al_1} \cdot \#\left([w] \cap
              \widetilde{\W}_{\alpha, \delta}\right).
    \end{equation*}
    Since $[w]$ is an orbit under $C_n$, an element $\tau \in C_n$
    fixes $w$ if and only if $\tau$ fixes $[w]$ pointwise. Thus, for any
    $\tau \in C_n$,
    \begin{equation}
      \label{eq:twiddle_fixed}
      \#\W_{\alpha, \delta}^\tau = \frac{n}{\alpha_1} \cdot
         \#\widetilde{\W}_{\alpha, \delta}^\tau.
    \end{equation}
    Combining \eqref{eq:twiddle_fixed} and
    Lemma~\ref{lem:fix_iff} now shows that for any $\tau \in C_g$,
    \begin{align}
      \label{eqn:SameNumFPs}
      \#\W_{\alpha, \delta}^\tau = \frac{n}{\alpha_1} \cdot \#
         \left(\prod_{\ell=2}^m \SSS_\ell \times \M_\ell
               \right)^\tau.
   \end{align}
   Hence, by \eqref{eqn:SameNumFPs}, the CSP in \eqref{eqn:CSPSMCW},
   and \eqref{CSP1},
     \[ \lp \W_{\alpha, \delta}, C_g,
         \frac{n}{\alpha_1} \widetilde{\W}_{\alpha, \delta}^{\maj}(q) \rp \]
   exhibits the CSP. By \eqref{eq:W_twiddle_W_maj}, $\frac{n}{\alpha_1}
   \widetilde{\W}_{\alpha, \delta}(q) \equiv \W_{\alpha, \delta}^{\maj}(q)$
    modulo $q^d - 1$, hence also modulo $q^g - 1$ since $ g \mid d $,
    completing the proof.
\end{pf}

We have now finished the verification of the conditions in
Lemma~\ref{lem:extend.csp} for $\W_{\alpha, \delta}$. Condition (i) is
Theorem~\ref{thm:CSPmodg}, Condition (ii) is Proposition~\ref{Propn:Periodg}, and
Condition (iii) is Lemma~\ref{lem:freqdivg}. This completes the proof of
Theorem~\ref{thm:alpha_delta}.

\section{Refinements of Binomial CSP's}
\label{sec:subsets}

A key step in the proof of Theorem~\ref{thm:CSPmodg} was
Theorem~\ref{thm:rsw_subs_msubs} due to Reiner--Stanton--White, which says that
the triples
    \[ \left(\binom{[0,n-1]}{k}, C_n, \binom{n}{k}_q\right) \qquad \text{and} \qquad
       \left(\mch{[0,n-1]}{k}, C_n, \mch{n}{k}_q\right)\]
exhibit the CSP. Indeed, \cite{MR2087303} contains two proofs, one via representation theory \cite[\S 3]{MR2087303} and another by direct calculation \cite[\S 4]{MR2087303}. In this section, we give two refinements of related CSP's involving an action of $ C_d $ on sets of subsets (Theorem~\ref{thm:Fixedgcds}) and multisubsets (Theorem~\ref{thm:FixedsizesCSPMult}) for all $ d \mid n $, using shifted sum statistics. Our proof of the subset refinement, Theorem~\ref{thm:Fixedgcds}, does not use Theorem~\ref{thm:rsw_subs_msubs}, so it can be used as an alternative proof of the subset case of Theorem~\ref{thm:rsw_subs_msubs}. Our method is inspired by the rotation of subintervals used by Wagon and Wilf in \cite[\S3]{MR1269164}.

\subsection{Cyclic Actions and Notation}
\label{subsec:actions}

We define two different cyclic actions of the cyclic group of order $ d $ on $ [0, n - 1] $ and induce these actions to $ \ch{[0,n - 1]}{k} $ and $ \mch{[0,n - 1]}{k} $. We also fix notation for the rest of the section.

\begin{nota}
  \label{not:subsets}
  Fix $ n \in \bZ_{\ge 1}, k \in \bZ_{\ge 0} $, and $ d \mid n $. Let
\[
    \SSS = \ch{[0,n - 1]}{k}, \qquad \M = \mch{[0,n - 1]}{k}
\]
For all $ j \in [1, \f{n}{d} ] $, let
    \[ I_d^j \coloneqq [(j - 1)d, jd - 1], \]
  which we call a \textit{$ d $-interval}. For any composition $ \al = (\al_1, \dots, \al_{n/d}) \vDash k $ with $ n/d $ parts, let
\begin{align}
\label{eq:Fixedsizes}
       \SSS_\al &\coloneqq \{ A \in \SSS : \#(A \cap I_d^j) = \al_j \; \tx{for all} \; j \}, \\
\label{eq:Fixedsizes2}
       \M_\al &\coloneqq \{ A \in \M : \#(A \cap I_d^j) = \al_j \; \tx{for all} \; j \},
\end{align}
where the intersection in \eqref{eq:Fixedsizes2} preserves the multiplicity of $A$. We also fix cyclic groups $ C_d $, $ C_d' $ of order $ d $ whose actions are described below.
\end{nota}

  Let $C_d$ act
  on $[0, n-1]$ by simultaneous rotation of $d$-intervals,
  which is generated by the permutation
  \begin{align}
    \si_d \coloneqq (0 \; 1 \; \dots \; (d - 1)) \dots
             ((n - d) \; (n - d + 1) \; \dots \; (n - 1))
  \end{align}
  in cycle notation. On the other hand, $C_n$ has a unique subgroup $ C_d' $ of
  order $d$ which also acts on $[0, n-1]$ and is generated by the
  permutation
  \begin{align}
    \si_n^{n/d}
    = \lp 0 \; \lp \f{n}{d} \rp
        \; \dots \lp n - \f{n}{d} \rp \rp
        \dots \lp \lp \f{n}{d} - 1 \rp
        \; \lp \f{2n}{d} - 1 \rp  \; \dots (n - 1) \rp.
  \end{align}
Induce these actions of $ C_d $ and $ C_d' $ up to $ \SSS $ and $ \M $ by
\begin{align*}
    g \dd \{ a_1, \dots, a_k \} \coloneqq \{ g \dd a_1, \dots, g \dd a_k \}.
\end{align*}
Notice that the action of $ C_d $ restricts to $ \SSS_\al $ and $ \M_\al $
for any $ \al = (\al_1, \dots, \al_{n/d}) \vDash k $.

Let $(G, X)$ be a pair where $G$ is a group acting on a set $X$.
A morphism of group actions $(G, X) \to (G', X')$ is a pair
$(\phi, \psi)$ where $\phi \colon G \to G'$ is a group homomorphism
and $\psi \colon X \to X'$ is a map of sets which satisfy
  \[ \psi(g \cdot x) = \phi(g) \cdot \psi(x)
     \; \tx{for all} \; g \in G, x \in X. \]

\begin{rem}
  \label{rem:Cdactions}
  The actions of $C_d$ and $C_d'$ on $[0, n-1]$ are isomorphic since
  $\sigma_d$ and $\sigma_n^{n/d}$ have the same
  cycle type. This isomorphism explicitly arises from
  $\phi \colon \sigma_d \mapsto \sigma_n^{n/d}$ with
  $\psi \colon 0 \mapsto 0$,
  $1 \mapsto \frac{n}{d}$, etc. Thus the actions of $C_d$ and $C_d'$
  on $ \SSS $ and $ \M $ are isomorphic as well.
\end{rem}

Recall the $ \Sum $ statistic sums the elements of a set or multiset.
We also use the following shifted sum statistic. For $A \in \SSS $, let
\begin{align}
  \label{eqn:shiftedsum}
  \Sum'(A) \coloneqq \sum_{a \in A} a - \sum_{i=0}^{k-1} i
            = \Sum(A) - \ch{k}{2}.
\end{align}
Recall from \eqref{eq:subset_sum} and \eqref{eq:multiset_sum} that
\begin{equation}
  \label{eq:sumprime}
  \SSS^{\Sum'}(q) = \binom{n}{k}_q,
  \qquad \M^{\Sum}(q) = \mch{n}{k}_q.
\end{equation}
Using \eqref{eq:sumprime}, we may restate Theorem~\ref{thm:rsw_subs_msubs}
as saying that
  \[ (\SSS, C_n, \SSS^{\Sum'}(q)) \qquad \text{and} \qquad
     (\M,C_n,\M^{\Sum}(q)) \]
exhibit the CSP. Moreover, under the restricted action of $C_d' \sub C_n$ on
$\M$ and $\SSS $,
  \[ (\SSS, C_d', \SSS^{\Sum'}(q)) \qquad \text{and} \qquad
     (\M, C_d', \M^{\Sum}(q)) \]
exhibit the CSP by Remark~\ref{rem:relatedCSPs}.
By Remark~\ref{rem:Cdactions},
\begin{equation}
  \label{eq:csp_sms_refinement}
  (\SSS, C_d, \SSS^{\Sum'}(q)) \qquad\text{and}\qquad
  (\M, C_d, \M^{\Sum}(q))
\end{equation}
also exhibit the CSP.

\begin{ex}
\label{Ex:Cdp1'action}
Let $ n = 8 $, $ k = 4 $, and $d = 4$. Abbreviating $\{0, 4, 5, 6\}$ as
$0456$, etc., gives
\begin{align*}
    \SSS_{(1,3)} = \{ & 0456,0457,0467,0567,
                      1456,1457,1467,1567,\\
                      & 2456,2457,2467,2567,
                      3456,3457,3467,3567
    \}.
\end{align*}
Here, $ C_4 $ acts on $[0, 8-1]$ by the permutation $ (0123)(4567) $, and $ C_4' $ acts by $ (0246)(1357) $. $ \M_{(1, 3)}$ contains $\SSS_{(1, 3)}$ in addition to, for instance, $0444$.
\end{ex}

\subsection{A Multisubset Refinement}
\label{subsec:multiset_refinement}

We next prove a refinement of the CSP triple $ (\M, C_d, \M^{\Sum}(q)) $ in \eqref{eq:csp_sms_refinement} by fixing sizes of intersections with the $ d $-intervals.

\begin{thm}
\label{thm:FixedsizesCSPMult}
Recall Notation~\ref{not:subsets}, and fix a composition $ \al = (\al_1, \dots, \al_{n/d}) \vDash k $. Then, $ (\M_\al, C_d, \M_\al^{\Sum}(q)) $ refines the CSP triple $ (\M, C_d, \M^{\Sum}(q)) $.
\end{thm}

\begin{pf} Separating the $ d $-intervals into different multisubsets gives
\begin{align}
    \M_\al \cong \mch{[0,d - 1]}{\al_1} \x \cdots \x \mch{[0,d - 1]}{\al_{n/d}},
\end{align}
which preserves the natural $ C_d $-action and $ \Sum $ statistic modulo $ d $. Since
\[
    \lp \mch{[0,d - 1]}{\al_j}, C_d, \mch{[0,d - 1]}{\al_j}^{\Sum}(q) \rp
\]
exhibits the CSP for all $ j $, the result follows from Remark~\ref{rem:relatedCSPs}.
\end{pf}

The following analogous result holds for subsets.

\begin{prop}
\label{prop:FixedsizesCSP}
Recall Notation~\ref{not:subsets}, and additionally fix a composition $ \al = (\al_1, \dots, \al_{n/d}) \vDash k $. Then $ (\SSS_\al, C_d, \SSS_\al^{\Sum^{\ast}}(q)) $ exhibits the CSP, where
\begin{align}
    \Sum^{\ast}(A) \coloneqq \Sum(A) - \sum_{ j = 1}^{n/d} \binom{\al_j}{2}.
\end{align}
\end{prop}

\begin{pf} Separating the $ d $-intervals into different subsets gives
\begin{align}
\label{eq:subsetprod}
    \SSS_\al \cong \ch{[0,d - 1]}{\al_1} \x \cdots \x \ch{[0,d - 1]}{\al_{n/d}},
\end{align}
which preserves the $ C_d $-action and $ \Sum $ statistic modulo $ d $. Since
\[
    \lp \binom{[0,d - 1]}{\al_j}, C_d, \binom{[0,d - 1]}{\al_j}^{\Sum - \ch{\al_j}{2}}(q) \rp
\]
exhibits the CSP for all $ j $, $ (\SSS_\al, C_d, \SSS_\al^{\Sum^{\ast}}(q)) $ exhibits the CSP by Remark~\ref{rem:relatedCSPs}.

\end{pf}

\begin{rem}
Since we must shift the $\Sum$ statistic by different amounts depending on $ \al $, Proposition~\ref{prop:FixedsizesCSP} is not a CSP refinement, in contrast to Theorem~\ref{thm:FixedsizesCSPMult}.
\end{rem}

\subsection{A Subset Refinement}
\label{subsec:subsets}

We next prove an honest refinement of the CSP triple
$ (\SSS, C_d, \SSS^{\Sum'}(q)) $ in \eqref{eq:csp_sms_refinement}. To do so, we
restrict to certain subsets of $S$ for each divisibility chain ending in
$n$. Our proof again inductively extends CSP's up from cyclic subgroups
of $C_d$ using Lemma~\ref{lem:extend.csp}. In this subsection we first define
our restricted subsets and give some examples. We then present a series
of lemmata verifying the conditions of Lemma~\ref{lem:extend.csp} before
proving our refinement, Theorem~\ref{thm:Fixedgcds}.

\begin{defn}
Suppose $ e \mid d \mid n $. Let
\begin{align}
     \G_{d, e} \coloneqq
        \{ A \in \SSS :  \gcd(d, \#(A \cap I_{d}^1), \#(A \cap I_{d}^2),
                                 \ldots, \#(A \cap I_{d}^{n/d}) ) = e \}.
\end{align}
\end{defn}

We have $ \G_{n, \gcd(n, k)} = \SSS$ and $ \G_{n, e} = \varnothing$ for all
other $e$.

\begin{rem}
  By conditioning on the sizes of the intersections with $ d $-intervals,
  $ \G_{d,e} $ decomposes as the disjoint union
  \begin{align}
    \label{eq:GasUnion}
    \G_{d,e} = \coprod
      \SSS_{\al},
  \end{align}
  ranging over all $ \al = (\al_1, \dots, \al_{n/d}) \vDash k$
  satisfying
  \begin{align*}
      \gcd(d, \al_1, \dots, \al_{n/d}) = e.
  \end{align*}
\end{rem}

\begin{ex}
  If $ n = 4, k = 2 $, then abbreviating $\{0, 2\}$ as $02$, etc., gives
  \begin{gather*}
    \G_{1,1} = \{01,02,03,12,13,23\} = \SSS, \\
    \G_{2,1} = \{02,03,12,13\}, \qquad
    \G_{2,2} = \{01,23\}, \\
    \G_{4,1} = \varnothing, \qquad
    \G_{4,2} = \{01,02,03,12,13,23\} = \SSS, \qquad
    \G_{4,4} = \varnothing.
  \end{gather*}
  Consequently,
  $\G_{4,2} \cap \G_{2,1} = \{02,03,12,13\}$ and
  $\G_{4,2} \cap \G_{2,2} = \{01,23\}$.
\end{ex}

\begin{defn}
  \label{def:G_D}
  Suppose $D$ is a totally ordered chain in the divisibility lattice
  ending with $\gcd(n, k) \mid n$,
  i.e.~$D = d_p \mid d_{p - 1} \mid \cdots \mid d_0 \mid n$ where
  $d_0 \coloneqq \gcd(n, k)$. Write
    \[ \G_D \coloneqq \G_{n,d_0} \cap \G_{d_0, d_1} \cap \cdots \cap
              \G_{d_{p - 1},d_p} \subset \SSS. \]
\end{defn}

We may now state our subset refinement. The proof is postponed to the
end of this subsection.

\begin{thm}
  \label{thm:Fixedgcds}
  Using Notation~\ref{not:subsets}, let $D$ be a totally ordered chain in the
  divisibility lattice ending with $\gcd(n, k) \mid n$ and starting with
  $e \mid d$. Then, $(\G_D, C_d, \G_D^{\Sum'}(q))$
  refines the CSP triple $(\SSS, C_d, \SSS^{\Sum'}(q))$.
\end{thm}

\begin{ex}
  \label{ex:Gsumdist}
  If $n = 4$, $k = 2$, and $D = 1 \mid 2 \mid 4$, then
  $ G = \G_{4,2} \cap \G_{2,1} $ has $ C_2 $ orbits $ \{ 02,13 \} $ and
  $ \{03,12\} $. Moreover,
  \[
     G^{\Sum'}(q) = q^1 + 2q^2 + q^3 \equiv 2(q^0 + q^1)
     \quad (\Mod q^2 - 1),
  \]
  so $ (G, C_2, G^{\Sum'}(q)) $ exhibits the CSP by \eqref{CSP2}.
\end{ex}

In fact, the subset case of Theorem~\ref{thm:rsw_subs_msubs} is the special case
$D = \gcd(n, k) \mid n$ of Theorem~\ref{thm:Fixedgcds}, so the proof below of
Theorem~\ref{thm:Fixedgcds} yields an alternative proof of the subset case of
Theorem~\ref{thm:rsw_subs_msubs}.

\begin{cor}
  $ (\SSS, C_n, \SSS^{\Sum'}(q)) $ exhibits the CSP.
\end{cor}

\begin{lem}
\label{lem:G_D_props}
  Let $D$ be a totally ordered chain in the divisibility lattice ending
  with $\gcd(n, k) \mid n$ and beginning with $e \mid d$. Suppose
  $C_e'$ is the unique subgroup of $C_d$ of order $e$.
  \begin{enumerate}[(i)]
    \item $ \G_D = \coprod \SSS_\alpha$, where the disjoint union is over
      a subset of the sequences $\alpha$
      satisfying $\alpha = (\alpha_1, \ldots, \alpha_{n/d}) \vDash k$
      and $\gcd(d, \alpha_1, \ldots, \alpha_{n/d}) = e$.
    \item $\G_D$ is closed under the $C_d$ and $C_e$-actions on $\SSS$.
    \item The $C_e'$ and $C_e$-actions on $\G_D$ are isomorphic.
    \item For any $C_d$-orbit $\cO$ of $\G_D$, we have
      $\frac{d}{|\cO|} \mid e$.
    \item The $\Sum'$ statistic has period $e$ modulo $d$ on $\G_D$.
  \end{enumerate}
\end{lem}

\begin{pf}
For (i), by \eqref{eq:GasUnion} we have $ \G_{d, e} = \coprod \SSS_\alpha$ where
  \[ \alpha = (\alpha_1, \ldots, \alpha_{n/d}) \vDash k \qquad \text{and} \qquad \gcd(d, \alpha_1, \ldots, \alpha_{n/d}) = e. \]
  Write $ D = e \mid d \mid d_1 \mid \dots \mid d_r = n $. For all $ j = 1, \dots, r $, since $ d \mid d_j $, each $ d_j $-interval is a union of $ d $-intervals. Thus, for $ A \in \SSS_\al $, whether $ A \in \G_{d_{j}, d_{j + 1}} $ for any $ j $ depends only on $ \al $, so $ \G_{d_{j}, d_{j + 1}} \cap \SSS_\al = \es $ or $ \G_{d_{j}, d_{j + 1}} \cap \SSS_\al = \SSS_\al $. Now (i) follows from $ \G_D = \G_{d,e} \cap \G_{d, d_1} \cap \dots \cap \G_{d_{r - 1}, d_r} $.

  For (ii), by (i) it suffices to show that each $ \SSS_\alpha$ is
  closed under the $C_d$ and $C_e$-actions. Since $\sigma_d$ rotates
  $d$-intervals, it preserves the size of each $d$-interval, so
  $\sigma_d$ indeed maps $ \SSS_\alpha$ to itself. The same argument applies
  with $\sigma_e$ in place of $\sigma_d$.

  For (iii), by (i), it suffices to show the $C_e$ and $C_e'$-actions
  on $ \SSS_\alpha$ are isomorphic. Recalling \eqref{eq:subsetprod}, we have
  \[
    \SSS_\al \cong \ch{[0,d - 1]}{\al_1} \x \cdots \x
    \ch{[0,d - 1]}{\al_{n/d}}.
  \]
  By Remark~\ref{rem:Cdactions}, the actions of $ C_e $ and $ C_e' $ on
  $ \ch{[0,d - 1]}{\al_j} $ are isomorphic for each $ j $, so their
  actions on $ \SSS_\al $ are isomorphic as well.

  For (iv), pick $A \in \cO$ with $A \in \SSS_\alpha$ for $\alpha$ as in
  (i). Let $A_j \coloneqq A \cap I_d^j$, which has $\alpha_j$
  elements. Viewing $A_j$ as a multiplicity word $ w_j $ as in
  Definition~\ref{defn:multword}, we see that $A_j$ has $d-\al_j$ zeros and $\al_j$
  ones. For all $ j $, $ w_j $ is some word repeated $ \f{d}{|\cO|} $
  times. Using the two-letter case of Lemma~\ref{lem:freqdivg}, we have
  $ \f{d}{|\cO|} \mid \freq(w_j) \mid \al_j $. Thus
  $ \f{d}{|\cO|} \mid \gcd(d, \al_1, \dots, \al_{n/d}) = e $.

  For (v), it suffices to show that $\Sum'$ has period $e$ modulo $d$
  on $\SSS_\alpha$ for $\alpha$ as in (i). By the gcd condition, there exist
  $ c_1, \ldots, c_{n/d} \in \bZ $ such that
  \[
     c_1 \al_1 + \dots + c_{n/d} \al_{n/d} \equiv e \; (\Mod d).
  \]
  For some particular $A \in \SSS_\al $, consider cyclically rotating the
  elements of $A \cap I_d^j$ forward by $c_j$ in $ I_d^j $ for all
  $ j $. The result is a bijection $\phi \colon \SSS_\al \to \SSS_\al $ that satisfies $ \Sum'(\phi(A)) \equiv \Sum'(A) + e \; (\Mod d) $, from which (v) follows.
\end{pf}

\begin{ex}
  Let $n=12, k=8$, and $D = 1 \mid 2 \mid 4 \mid 12$. Then
    \[ \G_D = \G_{12, 4} \cap \G_{4, 2} \cap \G_{2, 1}
           = \G_{4, 2} \cap \G_{2, 1}. \]
  We have $ \G_{2, 1} = \coprod \SSS_\alpha$ where $\alpha = (\alpha_1,
  \ldots, \alpha_6) \vDash 8$ and $\gcd(2, \alpha_1, \ldots, \alpha_6)
  = 1$. Similarly $ \G_{4, 2} = \coprod \SSS_\beta$ where $\beta = (\beta_1,
  \beta_2, \beta_3) \vDash 8$ and $\gcd(4, \beta_1, \beta_2, \beta_3)
  = 2$. In fact,
    \[ \varnothing \subsetneq \G_D \subsetneq \G_{2, 1} \]
  since, for instance, $\SSS_\alpha \subset \G_D$ when
  $\alpha = (4, 0, 1, 1, 1, 1)$ while $\SSS_\alpha \subset \G_{2, 1} - \G_D$
  when $\alpha = (2, 0, 2, 1, 2, 1)$.
\end{ex}

\begin{lem}
  \label{lem:G_dd_trivial}
  Let $d \mid n$. The $C_d$ action on $ \G_{d, d}$ is trivial and
    \[ (\G_{d, d}, C_d, \G_{d, d}^{\Sum'}(q)) \]
  exhibits the CSP.
\end{lem}

  \begin{pf}
    All subsets in $ \G_{d, d}$ have each $d$-interval either full or
    empty, so $ C_d $ fixes every $ A \in \G_{d,d} $. By \eqref{CSP2},
    $(\G_{d, d}, C_d, \G_{d, d}^{\Sum'}(q))$ thus exhibits the CSP if and
    only if $\G_{d, d}^{\Sum'}(q) \equiv \# \G_{d, d}$ mod $(q^d - 1)$.
    If $ \G_{d, d} = \varnothing$ the result is trivial, so take
    $ \G_{d, d} \neq \varnothing$. For any $A \in \G_{d, d}$, since each
    $d$-interval is full or empty, we have $d \mid k$ and
    \begin{align}
      \Sum'(A) \equiv
        \f{k}{d} \ch{d}{2}
        - \ch{k}{2} \equiv \f{k(d - k)}{2}
      \equiv 0 \; (\Mod d).
    \end{align}
  \end{pf}

We may now prove Theorem~\ref{thm:Fixedgcds}.

\begin{pf}[Proof of {Theorem~\ref{thm:Fixedgcds}}]
  We induct on $d$. If $d=1$, then the relevant triple $( \G_D, C_1, \G_D^{\Sum'}(q))$
  exhibits the CSP trivially. For the induction step, we first claim
  that $( \G_D, C_e, \G_D^{\Sum'}(q))$ exhibits the CSP. If $e=d$,
  then $ \G_D \sub \G_{d,d} $, so by Lemma~\ref{lem:G_dd_trivial} the $ C_e $
  action is trivial. It is easy to see that CSP's with trivial actions
  refine to arbitrary subsets, so $(\G_D, C_e, \G_D^{\Sum'}(q))$ exhibits
  the CSP in this case. If $e < d$, by
  conditioning on the sizes of the intersections of the $e$-intervals,
  we can write
  \begin{equation}
    \label{eqn:Nextgcd}
    \G_D = \coprod_{f \mid e} \G_{f \mid D}
  \end{equation}
  where $f \mid D$ denotes the chain with $f$ prepended to $D$. Hence
  $(\G_{f \mid D}, C_e, \G_{f \mid D}^{\Sum'}(q))$ exhibits the CSP by
  induction for each $f \mid e$, since $ f \mid D $ begins with
  $ f \mid e $. Thus $(\G_D, C_e, \G_D^{\Sum'}(q))$ exhibits the CSP by
  \eqref{eqn:Nextgcd}, proving the claim.

  In order to realize the $(\G_D, C_d, \G_D^{\Sum'}(q))$ CSP triple from the CSP triple $(\G_D, C_e, \G_D^{\Sum'}(q))$, we verify the conditions of
  Lemma~\ref{lem:extend.csp}. From Lemma~\ref{lem:G_D_props}(ii), the restriction
  of the $C_d$-action on $\G_D$ to the subgroup
  $C_e' \subset C_d$ of size $e$ is isomorphic to the $C_e$-action
  on $\G_D$, giving Condition (i). Condition (ii) is
  Lemma~\ref{lem:G_D_props}(v), and Condition (iii) is
  Lemma~\ref{lem:G_D_props}(iv). Thus $(\G_D, C_d, \G_D^{\Sum'}(q))$
  exhibits the CSP by Lemma~\ref{lem:extend.csp}.
\end{pf}

\section{The Flex Statistic}
\label{sec:flex}

We conclude by formalizing the notion of \textit{universal} sieving
statistics and giving an example, \textit{flex}, in the context of words.
We end with an open problem.

\begin{defn}
Given a set $W$ with a $C_n$-action, we say $\stat \colon W \to \bZ_{\geq 0}$ is a \textit{universal CSP statistic} for $ (W, C_n)$ if $ (\cO, C_n, \cO^{\stat}(q)) $ exhibits the CSP for all $C_n$-orbits $\cO$ of $W$.
\end{defn}

\begin{defn}
  Let $\lex(w)$  denote the index at which $w$ appears when
  lexicographically ordering the necklace $[w]$, starting from $0$. Let \textit{flex} be the product
    \[ \flex(w) \coloneqq \freq(w) \lex(w). \]
\end{defn}

For example, listing $ N = [15531553] $ in lexicographic order gives
\[
    N = \{15531553, 31553155, 53155315, 55315531 \},
\]
so, noting $ \freq(N) = 2 $, we have
\begin{align*}
    \lex(15531553) = 0, & \quad \flex(15531553) = 0, \\
    \lex(31553155) = 1, & \quad \flex(31553155) = 2, \\
    \lex(53155315) = 2, & \quad \flex(31553155) = 4,  \\
    \lex(55315531) = 3, & \quad \flex(55315531) = 6.
\end{align*}

\begin{lem}
  \label{lem:univCSP}
The function $\flex$ is a universal CSP statistic for $ (\W_n, C_n) $.
\end{lem}

\begin{pf} Let $ N $ be any necklace of length $ n $ words. Since $ \freq(N) = \f{n}{|N|} $, and $ \{ \lex(w) : w \in N \} = \{0, 1, \ldots, |N| -1\} $, we have
\begin{align}
\label{eqn:univCSP}
    N^{\flex}(q) = \sum_{j = 0}^{|N| - 1} q^{j \dd \f{n}{|N|}} = \f{q^{n} - 1}{q^{n/|N|} - 1},
\end{align}
so $ (N, C_n, N^{\flex}(q)) $ exhibits the CSP by \eqref{CSP2}.
\end{pf}

Given a universal sieving statistic $\stat$ on some set $W$, $\stat$ takes on precisely the values $\{0, n/d, \ldots, n-n/d\}$ modulo $n$ on any orbit of size $d$. The
converse holds as well. In this sense, up to shifting values by $n$,
universal sieving statistics are equivalent to total orderings on each
orbit $\cO$ of $W$.

Standing in contrast to Lemma~\ref{lem:univCSP}, $ (N, C_n, N^{\maj}(q)) $ does not exhibit the CSP when $N = [123123]$, so $\maj$ is not a universal CSP statistic on $ (\W_n, C_n) $. However, $\maj$ trivially refines to the orbit $N=\{1^n\}$ for any $n$. Since refinement is not generally closed under intersections, it is not clear if there is any useful sense in which $\maj$ on words can be ``maximally refined.''

It follows from Lemma~\ref{lem:univCSP} and \eqref{CSP2} that
Theorem~\ref{thm:alpha_delta} is equivalent to the following.

\begin{thm}
  \label{thm:majflex}
  The statistics $ \flex $ and $ \maj $ are equidistributed modulo $ n $
  on $ \W_{\al,\de} $.
\end{thm}

Indeed, we were originally led to Theorem~\ref{thm:alpha_delta}
through an exploration of the irreducible multiplicities of the
so-called higher Lie modules (see e.g.~\cite{MR1984625}), which uncovered
the fact that $ \flex $ and $ \maj $ are equidistributed modulo $ n $
  on $ \W_{\al} $. Data exploration led us to conjecture this equidistribution refined to fixed cyclic descent type as in Theorem~\ref{thm:majflex}. These connections will be explained
in a future publication. They also naturally suggest the problem of
finding explicit bijections proving Theorem~\ref{thm:majflex}, which we leave
as an open problem.

\begin{prob}
  For $\alpha \vDash n$ and $\delta$ any weak composition, find a
  bijection
    \[ \va \colon \W_{\al, \de} \to \W_{\al, \de} \]
  satisfying
  \begin{align}
    \label{eq:bij_maj_flex}
    \maj(\va(w)) \equiv \flex(w) \;(\Mod n).
  \end{align}
\end{prob}

\section*{Acknowledgments}

We were partially supported by the National Science Foundation grant DMS-1101017. We sincerely thank our advisor, Sara Billey, for her many helpful suggestions, including connections to cyclic sieving, and for her very careful reading of numerous versions of the manuscript. We also thank Vic Reiner for fruitful early discussions related to Lie multiplicities, Yuval Roichman and his collaborators for generously sharing their preprints, and the anonymous referees for their thoughtful comments.

\bibliography{refs}{}

\begin{thebibliography}{RSW04}

\bibitem[AER17]{1801.00044}
Ron~M. Adin, Sergi Elizalde, and Yuval Roichman.
\newblock Cyclic descents for near-hook and two-row shapes, 2017.

\bibitem[ARR15]{MR3281144}
Drew Armstrong, Victor Reiner, and Brendon Rhoades.
\newblock Parking spaces.
\newblock {\em Adv. Math.}, 269:647--706, 2015.

\bibitem[ARR17]{1710.06664}
Ron~M. Adin, Victor Reiner, and Yuval Roichman.
\newblock On cyclic descents for tableaux, 2017.

\bibitem[AS17]{as17}
Connor Ahlbach and Joshua~P. Swanson.
\newblock Refined cyclic sieving.
\newblock {\em S\'em. Lothar. Combin.}, 78B:Art. 48, 12, 2017.

\bibitem[BER11]{MR2837599}
Andrew Berget, Sen-Peng Eu, and Victor Reiner.
\newblock Constructions for cyclic sieving phenomena.
\newblock {\em SIAM J. Discrete Math.}, 25(3):1297--1314, 2011.

\bibitem[Cel98]{MR1637728}
Paola Cellini.
\newblock Cyclic {E}ulerian elements.
\newblock {\em European J. Combin.}, 19(5):545--552, 1998.

\bibitem[ER17]{MR3682729}
Sergi Elizalde and Yuval Roichman.
\newblock On rotated {S}chur-positive sets.
\newblock {\em J. Combin. Theory Ser. A}, 152:121--137, 2017.

\bibitem[Gup78]{MR495467}
Hansraj Gupta.
\newblock A new look at the permutations of the first {$n$} natural numbers.
\newblock {\em Indian J. Pure Appl. Math.}, 9(6):600--631, 1978.

\bibitem[Kly74]{klyachko74}
A.~A. Klyachko.
\newblock Lie elements in the tensor algebra.
\newblock {\em Siberian Mathematical Journal}, 15(6):914--920, 1974.

\bibitem[LP12]{1202.4015}
Thomas Lam and Alexander Postnikov.
\newblock Alcoved {P}olytopes {I}{I}, 2012.

\bibitem[Mac13]{MR1506186}
P.~A. MacMahon.
\newblock The {I}ndices of {P}ermutations and the {D}erivation {T}herefrom of
  {F}unctions of a {S}ingle {V}ariable {A}ssociated with the {P}ermutations of
  any {A}ssemblage of {O}bjects.
\newblock {\em Amer. J. Math.}, 35(3):281--322, 1913.

\bibitem[Mac95]{MR1354144}
I.~G. Macdonald.
\newblock {\em Symmetric functions and {H}all polynomials}.
\newblock Oxford Mathematical Monographs. The Clarendon Press, Oxford
  University Press, New York, second edition, 1995.
\newblock With contributions by A. Zelevinsky, Oxford Science Publications.

\bibitem[Pec14]{MR3207480}
Oliver Pechenik.
\newblock Cyclic sieving of increasing tableaux and small {S}chr\"oder paths.
\newblock {\em J. Combin. Theory Ser. A}, 125:357--378, 2014.

\bibitem[Pet05]{MR2181371}
T.~Kyle Petersen.
\newblock Cyclic descents and {$P$}-partitions.
\newblock {\em J. Algebraic Combin.}, 22(3):343--375, 2005.

\bibitem[PSV16]{MR3537922}
Timothy Pressey, Anna Stokke, and Terry Visentin.
\newblock Increasing tableaux, {N}arayana numbers and an instance of the cyclic
  sieving phenomenon.
\newblock {\em Ann. Comb.}, 20(3):609--621, 2016.

\bibitem[Rho10]{MR2557880}
Brendon Rhoades.
\newblock Cyclic sieving, promotion, and representation theory.
\newblock {\em J. Combin. Theory Ser. A}, 117(1):38--76, 2010.

\bibitem[RSW04]{MR2087303}
V.~Reiner, D.~Stanton, and D.~White.
\newblock The cyclic sieving phenomenon.
\newblock {\em J. Combin. Theory Ser. A}, 108(1):17--50, 2004.

\bibitem[Sag11]{MR2866734}
Bruce~E. Sagan.
\newblock The cyclic sieving phenomenon: a survey.
\newblock In {\em Surveys in combinatorics 2011}, volume 392 of {\em London
  Math. Soc. Lecture Note Ser.}, pages 183--233. Cambridge Univ. Press,
  Cambridge, 2011.

\bibitem[Sch03]{MR1984625}
Manfred Schocker.
\newblock Multiplicities of higher {L}ie characters.
\newblock {\em J. Aust. Math. Soc.}, 75(1):9--21, 2003.

\bibitem[Spr74]{MR0354894}
T.~A. Springer.
\newblock Regular elements of finite reflection groups.
\newblock {\em Invent. Math.}, 25:159--198, 1974.

\bibitem[Sta99]{MR1676282}
Richard~P. Stanley.
\newblock {\em Enumerative combinatorics. {V}olume 2}, volume~62 of {\em
  Cambridge Studies in Advanced Mathematics}.
\newblock Cambridge University Press, Cambridge, 1999.

\bibitem[Sta12]{MR2868112}
Richard~P. Stanley.
\newblock {\em Enumerative combinatorics. {V}olume 1}, volume~49 of {\em
  Cambridge Studies in Advanced Mathematics}.
\newblock Cambridge University Press, Cambridge, second edition, 2012.

\bibitem[WW94]{MR1269164}
Stan Wagon and Herbert~S. Wilf.
\newblock When are subset sums equidistributed modulo {$m$}?
\newblock {\em Electron. J. Combin.}, 1:Research Paper 3, approx.\ 15, 1994.

\end{thebibliography}
\bibliographystyle{alpha}

\end{document}